\documentclass[11pt,twoside]{article}

\usepackage[latin1]{inputenc}
\usepackage[T1]{fontenc}
\usepackage[french]{babel}
\usepackage{aeguill}
\usepackage{amssymb}

\usepackage{graphicx}
\usepackage[a4paper,pdftex=true]{geometry}
\usepackage{fancyvrb}

\usepackage{amssymb,array}
\usepackage{syntax}
\usepackage{mdwtab}
\usepackage{mathenv}
\usepackage{float}
\usepackage{lscape}
\usepackage{here}
\usepackage{psfig}
\usepackage[a4paper]{geometry}
\usepackage{amssymb}
\usepackage{bm}
\usepackage{fancyhdr}
\usepackage{graphicx}
\usepackage{natbib}

\title{\Large \bf \sc Maximum pseudo-likelihood estimator for nearest-neighbours Gibbs point processes
}

\author{\thispagestyle{empty} {\sc By Jean-Michel Billiot$^1$, Jean-Fran\c{c}ois Coeurjolly$^{1,2}$} \\ 
{\sc and Rémy Drouilhet$^1$} \\
{\it Labsad, University of Grenoble 2, France}}

\date{January 3, 2006}

\usepackage{natbib}
\input tcilatex

\newcommand{\VE}{V}

\newcommand{\Par}{\theta}                    % paramètre
\newcommand{\ParV}{\Vect{\theta}}            % paramètre vectoriel
\newcommand{\ParVi}[1]{\Vect{\theta}^{(#1)}}            % paramètre vectoriel
\newcommand{\Pari}[1]{{\theta}^{(#1)}}           
\newcommand{\ParVT}{\Vect{\theta}^\star}     % paramètre vectoriel Vrai
\newcommand{\SpPar}{\Vect{\Theta}}           % espace des paramètres
\newcommand{\Vois}{\mathcal{\Vect{V}}}

\newcommand{\PL}[3]{PL_{#1} \left( #2;#3 \right)}
\newcommand{\LPL}[3]{LPL_{#1} \left( #2;#3 \right)}

\newcommand{\Dom}[1]{\Lambda_{(#1)}}

\newcommand{\V}[1]{\VE\left(#1\right)}

\newcommand{\VI}[2]{\VE\left(#1|#2\right)}

\newcommand{\VPar}[2]{\VE \left(#1 ; #2  \right)}
\newcommand{\VIPar}[3]{\VE \left(#1|#2  ; #3 \right)}

\newcommand{\dVIPar}[4]{\frac{\partial \VE}{\partial \Par_{#1}} \left(#2|#3  ; #4 \right)}
\newcommand{\ddVIPar}[5]{\frac{\partial^2 \VE}{\partial \Par_{#1} \partial \Par_{#2}  } \left(#3|#4  ; #5 \right)}

\newcommand{\dLPL}[3]{ \mathbf{LPL}_{\Lambda_{#1}}^{(1)} \left( #2 ; #3  \right) }
\newcommand{\dLPLt}[2]{ \mathbf{LPL}_{\widetilde{\Lambda}}^{(1)} \left( #1 ; #2 \right) }

\newcommand{\Estn}[1]{ \widehat{#1}_n }

\newcommand{\Dt}{\widetilde{D}}

\newcommand{\Dp}[2]{\left( #1 | #2 \right) }

\usepackage{bm}
\usepackage{amssymb}
\newcommand{\Var}{
        \mathbb{V}\mbox{ar}
}
\newcommand{\Mat}[1]{
  \underline{\bm{#1}}
}

\newcommand{\tr}[1]{
 {#1}^{\!  T}
}

\newcommand{\Vect}[1]{
  \bm{#1}
}

\newcommand{\Esp}{
  \bm{E}
}

\def\RR{\hbox{I\kern-.1em\hbox{R}}}
\def\NN{\hbox{I\kern-.1em\hbox{N}}}
\def\ZZ{\hbox{Z\kern-.3em\hbox{Z}}}
\def\PP{\hbox{P\kern-.8em\hbox{P}}}

\newlength{\hauteur}
\newlength{\profond}

\renewcommand{\Citet}[1]{\citeauthor{#1} \citeyearpar{#1}}

\begin{document}

\bigskip

\maketitle

\footnotetext[1]{Supported by a grant from IMAG Project AlpB.  }
\footnotetext[2]{Corresponding author}

% 60K35 Special processes : interacting random processes,...
% 60G55 Stochastic processes : point processes
% 62M30 Inference from stochastic processes : spatial processes.

\newpage

\begin{center}
{\large \bf Abstract } 
\bigskip

{\small \begin{minipage}{12cm}
This paper is devoted to the estimation of a vector $\Vect{\theta}$ parametrizing an energy function associated to some ``Nearest-Neighbours'' Gibbs point process, via the pseudo-likelihood method. We present some convergence results concerning this estimator, that is strong consistency and asymptotic normality, when only a single realization is observed. Sufficient conditions are expressed in terms of the local energy function and are verified on some examples.
\end{minipage}}
\end{center}

\vspace*{.5cm}

\section{Introduction}

Gibbs point processes first appeared in the theory of statistical physics. Historical aspects of the mathematical theory are covered briefly in~\cite{Kallenberg83}. The importance of the Gibbs point process as a model building
principle became widely recognized through these works. Indeed, the
class of Gibbs point processes is interesting because it allows to introduce
and study interactions between points through the modelling of an associated
potential function. This resulting gain explains their use in statistical
physics~\cite{Ruelle69}, \cite{Feynman72} (when taking interactions between molecules in models of dilute
gases into account) or in ecology (when analysing competitions between plants). Within the mechanics statistics framework, Gibbs states are defined as solutions of the well known equilibrium equations refered to Dobrushin-Lanford-Ruelle (D.L.R.) equations~\cite{Dobrushin69}, \cite{Lanford69}.
One way to introduce Gibbs point processes consists in using a family 
of local specifications with respect to a weight process. The Preston's theorems~(\cite{Preston76}) used precisely this approach
in order to give sufficient conditions on local specifications for the existence of Gibbs states. 

%\bigskip

Many proposals tried to estimate the potential function from the available point pattern data generated by some Gibbs point processes. If the potential belongs to a parametric family model, the most well-known methodology is the use of the likelihood function. The main drawback of this approach is that the likelihood function contains an unknown scaling factor whose value depends on the parameters and which is difficult to calculate.
The first class of models on which the estimation of the maximum likelihood has been undertaken is the class of pairwise interaction point processes. \cite{Ogata84} developed the maximum likelihood estimation method based on numerical approximations of the likelihood. \cite{Penttinen84} used a similar approach while applying a Monte Carlo method in a way to solve the likelihood equation by the stochastic Newton-Raphson algorithm. \cite{Moyeed91} proposed another iterative procedure for estimating the maximum likelihood estimator. For maximum likelihood  by Markov chain Monte Carlo, see  \cite{Geyer92}, \cite{Geyer99} and for U.L.A.N. conditions for maximum likelihood estimator, see \cite{Mase92}. An alternative approach consists in avoiding to optimize the likelihood function (because of the scaling factor problem) and introducing a pseudo-likelihood function instead. This idea originated from \cite{Besag74} in the study of lattice processes. \cite{Besag82} further considered this method for pairwise interaction point process, while \cite{Jensen91} generalized it to the general class of Gibbs point processes, see \cite{Mase95}, \cite{Mase99}, \cite{Jensen94}, \cite{Guyon91} for asymptotic properties.
A third way is the Takacs-Fiksel estimation method (\cite{Takacs86} \cite{Fiksel88}), which relies on a characteristic property of Gibbs processes using Palm measure. Asymptotics properties of Takacs-Fiksel estimator are studied in \cite{Heinrich92}, \cite{Billiot97}. A comparison of these different procedures applied to the Strauss model is presented in \cite{Diggle94}. The non parametric setting has been undertaken by \cite{Glotzl81} and \cite{Diggle87} (and the references therein). \cite{Heikkinen99} proposed a semiparametric estimator based on Bayesian smoothing techniques. A general review of the problem of statistical inference on spatial point processes can be found in the recent monograph of \cite{Moller03}.

The present study is devoted to ``Nearest-Neighbour" Gibbs point models by combining stochastic geometry arguments~(\cite{Stoyan95}) and computational geometry ones~(\cite{Preparata88}, \cite{Edelsbrunner88}, \cite{Boissonnat94}). Such models are introduced by~\cite{Baddeley89} where the neighbourhood relation depends on the realization of the process. Sufficient conditions (expressed in terms of the energy function) for the existence of such processes are proposed in \cite{BBD3} and \cite{BBD4}, where some examples are also proposed. The main one is a pairwise interaction point process where the neighbourhood relation corresponds to the (slightly modified) Delaunay graph of the realization of the process. 

In this paper, we  study a pseudo-likelihood estimator for such processes. More precisely, our framework is restricted to stationary Gibbs point processes based on energy function related to some graph (for instance the Delaunay graph) such that the energy function is invariant by translation and such that the local energy function is stable and quasi-local (or local). The main results of this paper are convergence results (strong consistency and asymptotic normality) of maximum pseudo-likelihood estimators in this framework. These results are obtained  when only a single realization is observed. Sufficient conditions are expressed in terms of the local energy function (which makes the results quite general) for some large family of parametrized energy functions. Among the different parametrizations, the exponential family is considered. 

The paper is organized as follows. Section~\ref{sec-prel} is devoted to some background on Gibbs point processes and to the description of our framework. The statistical model and the pseudo-likelihood method are presented in Section~\ref{sec-statModel}. Consistency and asymptotic normality of the maximum pseudo-likelihood estimator are respectively proved in Section~\ref{sec-cons} and Section~\ref{sec-norm}. Finally, the different sufficient conditions ensuring convergence results are verified on some examples in Section~\ref{sec-ex}. A short simulation is presented to check the effectiveness of maximum pseudo-likelihood estimator.

\section{Background on Gibbs point processes}\label{sec-prel}

\subsection{Gibbs point processes}

We define ${\cal B}$, ${\cal B}_b$ to be respectively the Borel $\sigma $%
-field and the bounded Borel boolean ring.

Let $\Omega$ denotes the class of locally finite subsets of $\RR^d$.
In particular, an element $\varphi$ of $\Omega$, also called configuration (of points), could be represented as $\varphi =\sum_{i\in {\bf I\!N}}\delta _{x_i}$ which is a simple counting
Radon measure in $\RR^d$ (i.e. all the points $x_i$ of $\RR^d$ are
distinct) where for every $\Lambda \in {\cal B}\,,\,\delta
_{x}(\Lambda )=1_\Lambda(x)$ is the Dirac measure and $1_A(.)$ is the indicator function of a set $A$. 
This space $\Omega$ is equipped with the vague topology, that is to say the weak
topology for Radon measures with respect to the set of continuous functions
vanishing outside a compact set. We also define the $\sigma$-field ${\cal F}$ 
spanned by the maps $\varphi \longrightarrow \varphi (\Lambda )\,,\Lambda
\in {\cal B}_b$, where $\varphi(\Lambda)$ corresponds to the number of points of $\varphi$ in $\Lambda$ due to the Radon measure representation of $\varphi$. The set of all configurations in a measurable set $\Lambda\subset \RR^d$ will be denoted by $\Omega_\Lambda$ and the corresponding  $\sigma $-field
${\cal F}_\Lambda$ is similarly defined. Furthermore, for any $\Lambda\in{\cal B}_b$,
$$
\left(\Omega,{\cal F}\right)=\left(\Omega_\Lambda,{\cal F}_\Lambda\right)\times
\left(\Omega_{\Lambda^c},{\cal F}_{\Lambda^c}\right)$$  
where $\Lambda^c=\RR^d \setminus \Lambda$ denotes the complementary of $\Lambda$ in $\RR^d$.
Finally, $\Omega_f $ denotes the class of all finite  subsets of $\RR^d$.

A point process on $\RR^d$ is a $\Omega$-valued random variable, denoted by $\Phi$, with probability distribution $P$ on $(\Omega, {\cal F}).$ and the intensity measure $\Lambda_{p}$ of $ P$ is defined as a measure on ${\cal B}$ such that for any $D\in{\cal B}$
$$ \Lambda_{p}({D})\,=\,\displaystyle\int_{\Omega} \varphi({D})\,\,{P}(d\varphi).$$
In the stationary case, $\Lambda_{p}({D})\,=\, \lambda_{p}\nu({D})$
where the constant $\lambda_p$ is called the intensity of $P$ and $\nu$ is the Lebesgue measure on $\RR^{d}$.

A Gibbs point process is usually defined using a family
of local specifications with respect to a weight process (often a stationary
Poisson process with distribution ${Q}$ and intensity $\lambda _{Q}=1$). Let $\Lambda$ be a bounded region in $\RR^d$. For such a process, given some configuration $\varphi
_{\Lambda^c}$ on $\Lambda^c$, the conditional probability on 
$\Lambda $ is of the form, for any $Y\in {\cal F}$~: 
$$ {\Pi }_\Lambda (\varphi,Y)= \left\{ \frac 1{Z_\Lambda(\varphi )}\int_{\Omega_\Lambda}
\exp \left( -\VI{\psi }{ \varphi_{\Lambda^c}} \right) 1_Y(\psi
\cup\varphi _{\Lambda^c}){Q_\Lambda}(d\psi )\right\} 1_{{R}_\Lambda}(\varphi),
$$
where $$Z_\Lambda(\varphi ) =\int_{\Omega_\Lambda} \exp \left( -\VI{\psi}{\varphi _{\Lambda^c}} \right){Q_\Lambda}(d\psi )$$
is called the partition function and ${R}_\Lambda =\{\varphi \in \Omega \,~:\,0<Z_\Lambda
(\varphi )<\infty \}$.

Whereas the finite energy function $\V{\varphi}$ measures the cost of any configuration,
the local energy $\VI{\psi}{\varphi}$ is defined as the energy required to add the points of $\psi$ in $\varphi$:
$$
\VI{\psi}{\varphi} =\V{\psi\cup \varphi} - \V{\varphi}.
$$
Let us notice that when $\psi$ reduces to one point $x$, we denote by a slight abuse $\VI{x}{\varphi}$ instead of $\VI{\{x\}}{\varphi}$.
It is well known that the collection of probability kernels $(\Pi_\Lambda
)_{\Lambda \in {\cal B}_b}$ satisfies the set of compatibility and
measurability conditions which define a local specification in the Preston's sense (\cite{Preston76}). The main condition is the consistency~: 
$$
\Pi_\Lambda \Pi _{\Lambda {\bf ^{^{\prime }}}}=\Pi_\Lambda \quad \text{%
for}\quad \Lambda ^{^{\prime }}\subset \Lambda. 
$$

Notice that some conditions are needed to ensure the existence of a
probability measure ${P}$ with respect to any local energy $V$ and
any weight process that satisfies the so-called Dobrushin-Lanford-Ruelle
(D.L.R.) equations~: 
$$
{P}(Y|{\cal F}_{\Lambda ^c})(\varphi )=\Pi_\Lambda (\varphi
,Y)\quad \text{for}\,\,{P}\,\,\text{a.e. }\varphi \in \Omega \,\quad
\text{for any } \Lambda \in {\cal B}_b\text{ and } Y\in {\cal F}. 
$$
For the general theory of Gibbs point processes, the reader may refer to \cite{Kallenberg83,Daley88,Stoyan95} and the references therein.

\subsection{Campbell and Palm measures and Glötz Theorem}

The reduced Campbell measure ${\cal C}_{p}^{!}$ of $P$ is a measure on ${\cal B}\otimes{\cal F}$ such that for any ${D}\in{\cal B}$ and any $Y\in {\cal F}$
$${\cal C}_{p}^{!}({D}\times Y)\,=\,\displaystyle\int_{\Omega}\int_{D} 1_{Y}(\varphi -\delta_{x})\,\varphi(dx)\,\,P(d\varphi). $$
When some measurable function $h$ from $\RR^d \times \Omega$ on $\RR^d$ is given, the following equation is often called the refined Campbell theorem 
$$\displaystyle \int_{\Omega}\sum_{x\in \varphi}h(x,\varphi-\delta_{x}){P}(d\varphi)\quad = \quad \displaystyle \int_{\RR^d \times\Omega}h(x,\varphi){\cal C}_{p}^{!}(d(x,\varphi)).$$
If the intensity measure $\Lambda_{p}$ is $\sigma-$finite, then for $\Lambda_{p}$- a.a. $x\in\RR^d$, the distribution ${P}_{x}^{!}$ on $(\Omega,{\cal F})$ exists. It is unique for $\Lambda_{p}$- a.a. $x\in \RR^{d}$ and such that
$${\cal C}_{p}^{!}({D}\times Y)\,=\,\displaystyle\int_{D}{P}_{x}^{!}(Y)\,\,\Lambda_{p}(dx)\quad\mbox{for any}\quad {D}\in{\cal B}, \, Y\in{\cal F}. $$
Then ${P}_{x}^{!}$ is called the reduced Palm distribution of the point process $P$ with respect to point~$x$.  Intuitively, the Palm distribution ${P}_{x}$ is the conditional probability of configurations of the point process given that the point $x$ belongs to the realization $\varphi$. Therefore, we have 
$$\displaystyle \int_{\Omega}\sum_{x\in \varphi}h(x,\varphi-\delta_{x}){ P}(d\varphi)\quad = \quad \displaystyle \int_{\RR^d \times\Omega}h(x,\varphi){ P}_{x}^{!}(d\varphi)\Lambda_{p}(dx).$$
When the process is stationary, one may apply the previous equation by replacing the intensity measure $\Lambda_{p}$ by its expression in this case $\lambda_p\nu$, and in this framework, \cite{Glotzl80} proved that ${P}$ $\in 
{\cal G}_0\left( \VE\right)$ if and only if the reduced Campbell measure ${\cal C}_p^{!}$ is absolutely continuous with respect to 
$\nu \times {P}$ and~: 
$$
\frac{d{\cal C}_p^{!}}{d\left( \nu \times {P}\right) }\left(
x,\varphi \right) =\lambda _{{P}}\frac{d{P}_x^{!}}{d{P}}\left( \varphi \right) =\exp \left( -\VI{x}{\varphi}\right) 
$$
where $\lambda _{{P}}=\stackunder{\Omega }{\int }\exp \left(
-\VI x\varphi\right) {P}\left( d\varphi \right)$ is the intensity
of the process ${P}$. In the particular case when $\VI{x}\varphi =0$, the point process corresponds to the stationary Poisson process ${Q}$. We know from the Slivnyak's theorem that ${Q}_x^{!}={ Q}$ which is one way of characterizing such process.

\subsection{Description of some Gibbs models}

This paper is mainly  devoted to the statistical study of some nearest-neighbours Gibbs point processes first introduced in \cite{Baddeley89}. More precisely, we are interested in models based on energy function of the form
\begin{equation} \label{eq-delMod}
\V{\varphi} = \sum_{k=1}^3\sum_{\xi \in Del_k(\varphi)} u^{(k)}(\xi,\varphi),
\end{equation}
where $Del_k(\varphi)$ is the set of clique of order $k$ of the Delaunay graph defined just below. For some $\varphi \in \Omega$ in general position, one defines $Del_3(\varphi)$ by the unique decomposition into triangles $\psi$ in which the convex hull of the circle $C(\psi)$ does not contain any point of $\varphi\setminus\psi$. The Delaunay graph is then defined by the set of edges:
$$
Del_2(\varphi) = \cup_{\psi \in Del_3(\varphi)} \mathcal{P}_2(\psi).
$$
In order to ensure the existence of such Gibbs state in $\RR^d$, \cite{BBD3} prove that 
the local stability and quasilocality properties (only expressed in terms of the energy function) are sufficient conditions of Preston's Theorem.

Without any additional modification, the previous model does not satisfy the previous assumptions. We then introduce some subgraphs. First let us denote, for some triangle $\psi$, by $D(\psi)$ the diameter of the circle circumscribed of $\psi $ and by $\beta \left( \psi \right) $ the smallest angle of $\psi$. 

\begin{definition}
Given any $\beta _0\in ]0,\pi /3]$, we introduce the following particular subset of $Del_{3}\left( \varphi \right)$~:
$$
Del_{3,\beta }^{\beta _{0}}\left( \varphi \right) =\left\{ \psi\in
Del_{3}\left( \varphi \right) :\beta \left( \psi\right) >\beta _{0}\right\}.
$$
The $\beta$-Delaunay graph of order $\beta_0$ of  any configuration $\varphi$ is the Delaunay subgraph 
defined by~:
$$
Del_{2,\beta }^{\beta _{0}}\left( \varphi \right) =\dbigcup\limits_{\psi\in
Del_{3,\beta }^{\beta _{0}}\left( \varphi \right) }\mathcal{P}_2(\psi).
$$
\end{definition}
The model obtained by replacing the original Delaunay graph by the $\beta$-Delaunay subgraph of order $\beta_0$  in~(\ref{eq-delMod}) satisfies the previous sufficient conditions of Preston's Theorem. From now on, this model is called the $\beta$-Delaunay model. In this spirit, some other models may be defined (see {\it e.g.} \cite{BBD3}, \cite{BBD4}) but we advice the reader to keep in mind the $\beta$-Delaunay model as the main example in order to illustrate the statistical results developped in this work.

The framework of this paper is restricted to stationary Gibbs point processes based on energy function related to some graph, denoted $\mathcal{G}_2(\varphi)$ for some finite configuration $\varphi$ ($\mathcal{G}_k(\varphi)$ representing the set of cliques of order $k$), of the form
\begin{eqnarray} 
\V{\varphi}  &=&
\sum_{k=1}^{K_{max}} \left\{ \sum_{\xi \in G_k(\varphi ) } u^{(k)} (\xi ; \varphi ) 
\right\} \label{eq-V} \\
&=& \Pari{1}|\varphi| + \sum_{k=2}^{K_{max}} \left\{ \sum_{\xi \in G_k(\varphi ) } u^{(k)} (\xi ;\varphi ) \right\}, \quad \mbox{ when } u^{(1)}\equiv \Pari{1} \nonumber 
\end{eqnarray}
and satisfying Assumptions $\mathbf{E_{1}}$, $\mathbf{E_2^{loc}}$ or more generally $\mathbf{E_{2}^{qloc}}$, $\mathbf{E_{3}}$ defined by:
\bigskip \bigskip
\begin{itemize}
\item[$\mathbf{E_{1}}$] $\V{\cdot}$ is invariant by translation.
\item[$\mathbf{E_2^{loc}}$] Locality of the local energy: there exists some fixed range denoted by $D$ such that for any $\varphi \in \Omega$ one has
$$
\VI{0}{\varphi} = \VI{0}{\varphi \cap \mathcal{B}(0,D)}.$$
\item[$\mathbf{E_2^{qloc}}$] Quasi-locality of the local energy: there exists a nonnegative function $\varepsilon$ vanishing asymptotically such that for any $\varphi \in \Omega$ one has
$$
|\VI{0}{\varphi} - \VI{0}{\varphi \cap \mathcal{B}(0,D)}|<\varepsilon(D).$$
\item[$\mathbf{E_{3}}$] Stability of the local energy: there exists $K\geq 0$ such that for any $\varphi \in \Omega$, 
$$\VI{0}{\varphi}\geq -K.$$
\end{itemize}

This framework includes some classical point processes such as~: 
\begin{itemize}
\item models based on the usual complete graph $\mathcal{G}_2(\varphi)=\mathcal{P}_2(\varphi)$) with pairwise interaction function satisfying a hard-core or inhibition condition and with finite range.
\item $k$-nearest neighbours models with pairwise interaction function bounded and with finite range (see \cite{BBD2}).
\item Widom-Rowlinson or area interaction model.
\item \ldots
\end{itemize}

\section{Statistical model and inference method} \label{sec-statModel}

\subsection{Statistical model}

We consider Gibbs point processes with energy function $\VPar{\cdot}{\ParV}$ parametrized as follows 

As a statistical model we consider a parametrized version of~(\ref{eq-V}) where the different $u^{(k)}(\xi,\varphi)$ depend on a vector parameters $\ParVi{k}$ and then denoted from now by $u^{(k)}(\xi;\varphi,\ParVi{k})$. It is assumed that the vector of parameters $\ParV=\left(\theta_1,\ldots,\theta_{p+1}\right)=\left(\Pari{1},\ParVi{2},\ldots,\ParVi{K_{max}} \right)\in \SpPar$ where $\SpPar$ is an open bounded set of $\RR^{p+1}$. 

Our data consist in the realization of a point process with energy function $\VPar{\cdot}{\ParVT}$ in a domain $\Lambda\subset \RR^d$ satisfying Assumptions $\mathbf{E_1}$ to $\mathbf{E_3}$. Thus, $\ParVT$ is the true parameter to be estimated. The Gibbs measure will be denoted by $P_{\ParVT}$. From~(\ref{eq-V}), we obtain easily the energy to insert a point $x$ in a configuration $\varphi$.

\begin{equation} \label{eq-VITheta}
\VIPar{x}{\varphi}{\ParV}  = \sum_{k=1}^{K_{max}} \left\{ \sum_{\xi \in G_k(\varphi \cup \{x \}) \setminus G_k(\varphi)} u^{(k)} (\xi ; \ParVi{k}, \varphi \cup \{x\}) -
\sum_{\xi \in G_k(\varphi) \setminus G_k(\varphi \cup \{ x \})} u^{(k)} (\xi ; \ParVi{k}, \varphi )
\right\}
\end{equation}

Theoretical results presented in the next sections are valid for a general energy function $\VPar{\cdot}{\ParV}$. But among this class of models, we will focus on energy functions described by~(\ref{eq-V}) and such that
$$u^{(k)}(\xi;\ParVi{k},\varphi)= \tr{\ParVi{k}}\Vect{u}^{(k)}(\xi;\varphi).$$
The energy can be rewritten
\begin{equation} \label{eq-modExp}
\VPar{\varphi}{\ParV} = \sum_{k=1}^{K_{max}} \sum_{\xi \in G_k(\varphi)} \tr{\ParVi{k}} \Vect{u}^{(k)} (\xi;\varphi) = \sum_{k=1}^{K_{max}}\tr{\ParVi{k}} \Vect{u}^{(k)}(\varphi) 
\end{equation}
where for any finite configuration $\varphi$ 
$$
\Vect{u}^{(k)}(\varphi) = \sum_{\xi \in G_k(\varphi)} \Vect{u}^{(k)} (\xi;\varphi) 
\quad \mbox{ and } \quad
\Vect{u}(\varphi) =\left(u_1(\varphi),\ldots,u_{p+1}(\varphi)\right)= \left( \Vect{u}^{(1)}(\varphi),\ldots, \Vect{u}^{(K_{max})}(\varphi)\right)
$$
For two finite configurations $\varphi$ and $\psi$, by denoting 
 \begin{equation} \label{eq-uPhiPsi}
\Vect{u}^{(k)}(\psi | \varphi) = \Vect{u}^{(k)}(\psi \cup \varphi) - \Vect{u}^{(k)}(\varphi )
\quad \mbox{ and } \quad
\Vect{u}(\psi|\varphi ) = \left( \Vect{u}^{(1)}(\psi|\varphi),\ldots, \Vect{u}^{(K_{max})}(\psi|\varphi)\right)
\end{equation}
we have for any point $x$
\begin{equation} \label{eq-Vexp}
\VPar{\varphi}{\ParV} = \sum_{k=1}^{K_{max}} \tr{\ParVi{k}} \Vect{u}^{(k)}(\varphi) =
\tr{\ParV} \Vect{u}(\varphi) 
\;\;\mbox{ and } \;\;
\VIPar{x}{\varphi}{\ParV} = \sum_{k=1}^{K_{max}} \tr{\ParVi{k}} \Vect{u}^{(k)}(x | \varphi) = \tr{\ParV} \Vect{u}(x|\varphi),
\end{equation}
where $\Vect{u}(x|\varphi)={(u_1(x|\varphi,\ldots,u_{p+1}(x|\varphi)}={(\Vect{u}^{(1)}(x|\varphi),\ldots,\Vect{u}^{K_{max}}(x|\varphi))}$. The local specification of the Gibbs point process associated to an energy function defined by~(\ref{eq-Vexp}) belongs to an exponential family.

\subsection{Pseudo-likelihood}

As precised in the introduction, the idea of maximum pseudo-likelihood is due to \cite{Besag75} who first introduced the concept for Markov random fields in order to avoid the normalizing constant. This work was then widely extended and \cite{Jensen91} (Theorem 2.2) obtained a general expression for Gibbs point processes. With our notation and up to a scalar factor the pseudo-likelihood defined for a configuration $\varphi$ and a domain of observation $\Lambda$ is denoted by $\PL{\Lambda}{\varphi}{\ParV}$ and given by
\begin{equation} \label{defPLpap}
\PL{\Lambda}{\varphi}{\ParV} = \exp \left( - \int_{\Lambda} \exp\left( - \VIPar{x}{\varphi}{\ParV} \right) dx \right)
  \prod_{x\in \varphi_\Lambda} \exp\left( -\VIPar{x}{\varphi \setminus x}{\ParV} \right).
\end{equation}
It is more convenient to define (and work with) the log-pseudo-likelihood function, denoted by $\LPL{\Lambda}{\varphi}{\ParV}$.
\begin{equation} \label{logPL}
\LPL{\Lambda}{\varphi}{\ParV}=
- \int_{\Lambda} \exp\left( - \VIPar{x}{\varphi}{\ParV} \right) dx  - \sum_{x\in\varphi_\Lambda}  \VIPar{x}{\varphi \setminus x}{\ParV} 
\end{equation}

\subsection{Main statistical tools}

Let us start by presenting a particular case of Campbell Theorem combined with Glötz Theorem that is widely used in our future proofs. For some finite configuration $\varphi$ (resp. for some set $G$) and for all $x$, we denote by $\varphi_x$ (resp. $G_x$) the configuration $\varphi$ (resp. the set $G$) translated of $x$.

\begin{corollary} \label{cor-CampGlotz}
If the probability measure $P$ is stationary and if the function $h(\cdot,\cdot)$ (used in Campbell Theorem) can be decomposed into $h(x,\varphi)=\mathbf{1}(x\in \Lambda) g(x,\varphi)$ for $\Lambda\subset \RR^d$ where $g(\cdot,\cdot)$ is such that $g(x,\varphi_x)=g(0,\varphi)$ for all $x$, then the refined Campbell theorem combined with Glötz Theorem allow us to obtain
\begin{equation} \label{eq-CamGlotz}
\Esp_P \Big( \sum_{x \in \Phi_\Lambda\setminus x} g(x,\Phi \setminus x) \Big) = |\Lambda| \; \Esp_P\Big(\; g(0,\Phi) \exp(-\VI{0}{\Phi})\; \Big)
\end{equation}
\end{corollary}

Let us now present a version of an ergodic theorem obtained by \cite{Nguyen79} and widely used in this paper. Let $\widetilde{D}>0$ and denote by $\Lambda_0$ the following fixed domain
$$
\Lambda_0 = \left\{ 
z \in \RR^2, -\frac{\widetilde{D}}2 \leq |z| \leq \frac{\widetilde{D}}2 
\right\},
$$
where for all $z\in \RR^2$, $|z|=\max(z_1,z_2)$.

\begin{theorem}[\Citet{Nguyen79}] \label{thmNguyen}
Let $\{H_G, G\in \mathbf{\mathcal{B}}_b\}$ be a family of random variables, which is covariant, that for all $x \in \RR^d$, 
$$
H_{G_x}(\varphi_x) = H_G(\varphi), \;\; a.s.
$$
and additive, that is for every disjoint $G_1,G_2 \in \mathbf{\mathcal{B}}_b$,
$$
H_{G_1\cup G_2} = H_{G_1} + H_{G_2} , \quad a.s.
$$
Let $\mathcal{I}$ be the sub$-\sigma-$algebra of $\mathcal{F}$ consisting of translation invariant (with probability 1) sets. Assume there exists a nonnegative and integrable random variable $Y$ such that $|H_G|\leq Y$ a.s. for every convex $G\subset \Lambda_0$. Then,
$$
\lim_{n \to +\infty} \frac1{|G_n|} H_{G_n} = \frac1{|\Lambda_0|}E (H_{\Lambda_0}|\mathbf{\mathcal{I}}) , \quad a.s.
$$
for each regular sequence $G_n \to \RR^d$.
\end{theorem}

\section{Consistency of the maximum pseudo-likelihood estimator} \label{sec-cons}

Maximizing the pseudo-likelihood is equivalent to minimize $U_n(\ParV)$ defined by
$$
U_n(\ParV) = - \frac1{|\Lambda_n|} \LPL{\Lambda_n}{\varphi}{\ParV}.
$$
We denote by $\Estn{ \ParV }=\Estn{\ParV}(\varphi)$ the maximum pseudo-likelihood estimator based on the configuration $\varphi$, alternatively defined as
$$
\Estn{\ParV}(\varphi) =  \mbox{argmin}_{\ParV \in \SpPar} U_n(\ParV)
$$

In this section, the existence of an ergodic measure is ensured, relatively to our framework, by Assumptions $\mathbf{E_1}$, $\mathbf{E_2^{qloc}}$ and $\mathbf{E_3}$.
The following Assumptions are needed to derive the almost sure convergence of this estimator.\\ 
\begin{itemize}
\item[$\mathbf{C_1}$] $(\Lambda_n)_{n\geq 1}$ is a regular sequence of domains such that $\Lambda_n \to \RR^2$ as $n \to +\infty$.
\item[$\mathbf{C_2}$] For all $\ParV \in \SpPar$, 
$$\VIPar{0}{\cdot}{\ParV} \in L^1(P_{\ParVT}).$$
\item[$\mathbf{C_3}$ ] For all $\ParV \in \Vect{\SpPar}\setminus\ParVT$
$$
P_{\ParVT} \Big( \left\{
\varphi, \;\; \VIPar{0}{\varphi}{\ParV} \neq \VIPar{0}{\varphi}{\ParVT}
\right\}\Big)>0$$
\item[$\mathbf{C_4}$] For all $\ParV,\ParV^\prime \in \SpPar$, there exists $c>0$ such that $P_{\ParVT}-$almost surely, we have
\begin{equation} \label{VThetaLip}
| \VIPar{0}{\Phi}{\ParV} - \VIPar{0}{\Phi}{\ParV^\prime} | \leq ||\ParV - \ParV^\prime||^c g(0,\Phi) 
\end{equation}
where $g(\cdot,\cdot)$ is a function such that for all $x$, $g(0,\Phi)=g(x,\Phi_x)$ and such that  $g(0,\cdot) \in L^1(P_{\ParVT})$.
\end{itemize}

\begin{remark} \label{rem-CsansLS} If one only assumes the existence of an ergodic measure, in particular without Assumption $\mathbf{E_3}$ (taken into account to express Assumptions $\mathbf{C_2}$ and $\mathbf{C_4}$) then 
\begin{list}{$\bullet$}{}
\item the condition $\mathbf{C_2}$ becomes~: for all $\ParV \in \Vect{\Theta}$, the variables $\VIPar{0}{\cdot}{\ParV}\exp\big(- \VIPar{0}{\cdot}{\ParVT}\big)$ and $\exp\big(-\VIPar{0}{\cdot}{\ParV} \big)$ are $P_{\ParVT}$-integrable.
\item the function $g(\cdot,\cdot)$ occuring in Assumption $\mathbf{C_4}$ is now such that for all $\ParV \in \Vect{\Theta}$, 
$g(0,\cdot)\exp\big(- \VIPar{0}{\cdot}{\ParV}\big) \in L^1(P_{\ParVT})$.
\end{list}
These Assumptions have been verified in \cite{Mase95} for the Ruelle class of pairwise interaction function with $\Vect{\theta}=(\beta,z)$ where $\beta$ represents the inverse temperature and $z$ the chemical potential.
\end{remark}

\begin{proposition} \label{prop-convMPLE}
Assume $P_{\ParVT}$ stationary, then under Assumptions $\mathbf{C_1}$ to $\mathbf{C_4}$, we have $P_{\ParVT}-$almost surely, as $n \to +\infty$
\begin{equation} \label{eq-convMPLE}
\Estn{\ParV}(\Phi) \rightarrow \ParVT
\end{equation}
\end{proposition}

Due to the decomposition of stationary measures as a mixture of ergodic measures (see \cite{Preston76}), one only needs to prove Proposition~\ref{prop-convMPLE} by assuming that $P_{\ParVT}$ is ergodic. Therefore, in Lemmas~\ref{lem-convUn} to \ref{lem-modCont}, $P_{\ParVT}$ is assumed to be ergodic.

The tool used to obtain the almost sure convergence is a convergence theorem for minimum contrast estimators established by \Citet{Guyon92}. Define
$$ K_n(\ParV,\ParVT)  = U_n(\ParV)- U_n(\ParVT)$$

\begin{lemma} \label{lem-convUn}
For all $\ParV\in \SpPar$, under Assumptions $\mathbf{C_1}$ and $\mathbf{C_2}$, we have $P_{\ParVT}-$almost surely, as $n \to +\infty$ 
\begin{equation} \label{defU}
U_n(\ParV) \rightarrow U(\ParV) =  \Esp_{P_{\ParVT} } \Big(  \exp\left(
- \VIPar{0}{\Phi}{\ParV} \right) \; + \;
\VIPar{0}{\Phi}{\ParV} \exp\big( - \VIPar{0}{\Phi}{\ParVT}\big)
\Big)
\end{equation}
\end{lemma}

\begin{proof}
Under Assumptions $\mathbf{C_1}$, $\mathbf{C_2}$ and $\mathbf{E_1}$, one can apply Theorem~\ref{thmNguyen} (\cite{Nguyen79}) to the process 
$$
H_{1,\Lambda_n}= \int_{\Lambda_n} \exp\left( -\VIPar{x}{\varphi}{\ParV} \right) dx .
$$
And from Corollary~\ref{cor-CampGlotz}, we obtain $P_{\ParVT}-$almost surely as $n \to +\infty$
\begin{equation} \label{H1Lambda}
  \frac1{|\Lambda_n|} H_{1,\Lambda_n}
 \rightarrow  \Esp_{P_{\ParVT} }\exp\left(- \VIPar{0}{\Phi}{\ParV} \right). 
\end{equation}
Now, define 
$$
H_{2,\Lambda_n} = \sum_{x \in \Phi_{\Lambda_n}} \VIPar{x}{\Phi\setminus x}{\ParV}
$$
Let $G \subset \Lambda_0$, we clearly have
$$
|H_{2,G} |\leq \sum_{x \in \Phi_G} | \VIPar{x}{\Phi\setminus x}{\ParV} | \leq
\sum_{x \in \Phi_{\Lambda_0}} | \VIPar{x}{\Phi\setminus x}{\ParV} |
$$
Under Assumption $\mathbf{C_2}$ and from Corollary~\ref{cor-CampGlotz}, we have 
$$\Esp_{P_{\ParVT}} \left( \sum_{x \in \Phi_{\Lambda_0}} | \VIPar{x}{\Phi\setminus x}{\ParV} | \right) = |\Lambda_0| \Esp_{P_{\ParVT}} \left( | \VIPar{0}{\Phi}{\ParV} | \exp\Big( - \VIPar{0}{\Phi}{\ParVT} \Big)\right) <+\infty
$$
This means that for all $G \subset \Lambda_0$, there exists a random variable $Y \in L^1(P_{\ParVT})$ such that $|H_{2,G}|\leq Y$. Thus, under Assumption $\mathbf{C_1}$ and from Theorem~\ref{thmNguyen} (\cite{Nguyen79}) and from Corollary~\ref{cor-CampGlotz}, we have $P_{\ParVT}-$almost surely
\begin{equation} \label{H2Lambda}
 \frac1{|\Lambda_n|} H_{1,\Lambda_n} \rightarrow  
\frac1{|\Lambda_0|} \Esp_{P_{\ParVT} }  \Big( \sum_{x \in \Phi_{\Lambda_0}} \VIPar{x}{\Phi\setminus x}{\ParV} \Big) = 
\Esp_{P_{\ParVT} }  \left(\VIPar{0}{\Phi}{\ParV} \exp\Big( - \VIPar{0}{\Phi}{\ParVT} \Big) \right).
\end{equation}
We have the result by combining (\ref{H1Lambda}) and (\ref{H2Lambda}).
\end{proof}

\begin{lemma} \label{lem-contraste}
Under the conditions of Lemma~\ref{lem-convUn}, the function $U_n(\cdot)$ defines a constrast function, that is there exists a function $K(\cdot,\ParVT)$ such that $P_{\ParVT}-$almost surely the following holds for all $\ParV \in \SpPar$:
$$
K_n(\ParV,\ParVT) \rightarrow K(\ParV,\ParVT)
$$
where $K(\cdot,\ParVT)$ is a positive which, under Assumption $\mathbf{C_3}$ is zero if and only if $\ParV=\ParVT$.
\end{lemma}

\begin{proof}
From Lemma~\ref{lem-convUn}, the function $K(\ParV,\ParVT)\geq 0$ can be written 
\begin{equation} \label{contrasteK}
K(\ParV, \ParVT)= \Esp_{P_{\ParVT}} \Big(
\exp(-\VIPar{0}{\Phi}{\ParVT}) \Big( 
\exp(\VIPar{0}{\Phi}{\ParV}-\VIPar{0}{\Phi}{\ParVT}) - (1+\VIPar{0}{\Phi}{\ParV}-\VIPar{0}{\Phi}{\ParVT})
\Big) \Big)
\end{equation}
The result is obtained using Assumption $\mathbf{C_3}$ and by noting that the function 
$t \mapsto \exp(t) -(1+t)$ is positive and is zero if and only if $t=0$.
\end{proof}

\begin{lemma} \label{lem-modCont}
Under Assumption $\mathbf{C_4}$, the functions $\ParV \mapsto U_n(\ParV)$ and $\ParV \mapsto K(\ParV, \ParVT)$ are continuous in $\ParV$. Moreover, the modulus of continuity of $U_n(\ParV)$ defined by
$$
W_n(\eta) = \sup \left\{ 
\Big|U_n(\ParV) - U_n(\ParV^\prime) \Big|, \ParV,\ParV^\prime \in \SpPar, || \ParV - \ParV^\prime || \leq \eta
\right\}
$$
is such that there exists a sequence $(\varepsilon_k)_{k \geq 1}$, with $\varepsilon_k \to 0$
as $k \to +\infty$ such that for all $k \geq 1$
\begin{equation} \label{modCont}
P_{\ParVT}\left( \limsup_{n \to +\infty}  \left( 
W_n \left(\frac1k\right) \geq \varepsilon_k
\right)\right) = 0.
\end{equation}
\end{lemma}

\begin{proof}
Under Assumption $\mathbf{C_4}$, it is sufficient to prove~(\ref{modCont}). Denote by 
$$
W_{1,n}\Big(\frac1k \Big) = \sup \left\{ 
\Big|\frac1{|\Lambda_n|} \int_{\Lambda_n}\!\! \Big( 
\exp\Big( -\VIPar{x}{\Phi}{\ParV} \Big) -\exp\Big( -\VIPar{x}{\Phi}{\ParV^\prime} \Big)
\Big) dx\Big| 
, \ParV,\ParV^\prime \in \SpPar, || \ParV - \ParV^\prime || \leq \frac1k
\right\}
$$
and
$$
W_{2,n}\Big(\frac1k \Big) = \sup \left\{ 
\Big|
\sum_{x \in \Phi_{\Lambda_n}} 
\VIPar{x}{\Phi\setminus x}{\ParV}-
\VIPar{x}{\Phi \setminus x}{\ParV^\prime}
\Big|, \ParV,\ParV^\prime \in \SpPar, || \ParV - \ParV^\prime || \leq \frac1k
\right\}
$$
Under Assumptions $\mathbf{E_3}$ and $\mathbf{C_4}$, one can prove that $P_{\ParVT}-$almost surely
\begin{eqnarray}
W_{1,n} \left( \frac1k \right)&\leq& \frac{\exp(K)}{k^c} \frac{1}{|\Lambda_n|} \int_{\Lambda_n} g(x,\Phi) dx \nonumber \\
W_{2,n} \left( \frac1k \right)& \leq & \frac{1}{k^c} \frac{1}{|\Lambda_n|} \sum_{x \in \Phi_{\Lambda_n}} g(x,\Phi \setminus x) \nonumber
\end{eqnarray}
Since $g(0,\cdot) \in L^1(P_{\ParVT})$, from Theorem~1 (\cite{Nguyen79}) there exists $N_0\in \NN$ such that for all $n \geq N_0$ we have 
$$W_{1,n} \left( \frac1k \right)\leq \frac{2 \exp(K)}{k^c} \Esp_{ P_{\ParVT}} (g(0,\Phi)) \quad \mbox{ and } \quad
W_{2,n} \left( \frac1k \right) \leq  \frac{2\exp(K)}{k^c}  \Esp_{ P_{\ParVT}} (g(0,\Phi)).$$
And so for all $n \geq N_0$, 
$$
W_n\left( \frac1k \right) \leq  \frac{\delta}{k^c} \qquad \mbox{ with } \delta= 4 \exp(K)\Esp_{ P_{\ParVT}} (g(0,\Phi)).
$$
Since
$$
\limsup_{n\to +\infty}  \left\{W_n\left( \frac1k \right) \geq \varepsilon_k \right\} \subset \left\{  \frac{\delta}{k^c} \geq \varepsilon_k \right\}.
$$
Thus, it is sufficient to choose $\varepsilon_k=\delta^\prime k^{-c}$ with $\delta^\prime>\delta$ to obtain the result.
\end{proof}

\noindent {\bf Proof of Proposition~\ref{prop-convMPLE}}

Lemma~\ref{lem-contraste} and Lemma~\ref{lem-modCont} ensure the fact that we can apply  Property 3.6 of \Citet{Guyon92} which asserts almost sure convergence for minimum contrast estimators. $\blacksquare$

The following proposition describes conditions $\mathbf{C_2},\mathbf{C_3}$ and $\mathbf{C_4}$ in the case of an exponential family. New conditions are denoted $\mathbf{C_2^{exp}}$, $\mathbf{C_3^{exp}}$ and $\mathbf{C_4^{exp}}$. For this result, let us consider energy functions described by (\ref{eq-Vexp}).

\begin{proposition} \label{prop-C234exp} ${ }$
Conditions $\mathbf{C_2}$ and $\mathbf{C_4}$ (resp. $\mathbf{C_3}$) can be replaced by $\mathbf{C_{2,4}^{exp}}$ (resp. $\mathbf{C_3^{exp}}$) where
\begin{itemize}
\item[$\mathbf{C_{2,4}^{exp}}$] There exists $\varepsilon>0$ such that for all $i=1,\ldots,p+1$
$$ u_i(0|\varphi) \in L^{1+\varepsilon}(P_{\ParVT}).$$
\item[$\mathbf{C_3^{exp}}$] Identifiability condition~: There exists $A_1,\ldots,A_{p+1}$, $p+1$ disjoint events of $\Omega$ such that $P_{\ParVT}(A_i)>0$ and such that for all $\varphi_1,\ldots,\varphi_{p+1} \in A_1 \times \cdots \times A_{p+1}$ the $(p+1)\times (p+1)$ matrix with entries $u_j(0|\varphi_i)$ is constant and invertible.
\end{itemize}
\end{proposition}

\begin{proof}
\begin{list}{$\bullet$}{}
\item Denote by $||\cdot||_{q}$ the norm defined for $z\in \RR^p$ by $||z||_q=\left(\sum_{i=1}^p |z_i|^q \right)^{1/q}$ with the obvious notation $||\cdot||=||\cdot||_2$. We have from Hölder's inequality
$$
| \VIPar{0}{\Phi}{\ParV} - \VIPar{0}{\Phi}{\ParV^\prime} | = \left| \left( \ParV-\ParV^\prime \right) \Vect{u}(0|\Phi) \right|\leq ||\ParV-\ParV^\prime||_{\frac{1+\varepsilon}{\varepsilon}} ||\Vect{u}(0|\Phi)||_{1+\varepsilon}.
$$ 
Since, $\SpPar$ is a bounded set there exists a constant $\kappa=\kappa(\varepsilon,\SpPar)$ such that we have 
$$
||\ParV-\ParV^\prime||_{\frac{1+\varepsilon}{\varepsilon}} \leq ||\ParV-\ParV^\prime||_2^{\frac{\varepsilon}{1+\varepsilon}} \; \;  ||\ParV-\ParV^\prime||_{\frac2{\varepsilon}}^{\frac{1}{1+\varepsilon}} \leq \kappa \; ||\ParV-\ParV^\prime||_2^{\frac{\varepsilon}{1+\varepsilon}} .
$$
Thus, we have (\ref{VThetaLip}), with $c=\frac{\varepsilon}{1+\varepsilon}$ and $g(0,\cdot)=||\Vect{u}(0|\Phi)||_{1+\varepsilon}$. And so, $\mathbf{C_{2,4}^{exp}}$ implies $\mathbf{C_4}$ and obviously $\mathbf{C_2}$.
\item Assumption $\mathbf{C_3^{exp}}$ means that for all $\Vect{y}\in \RR^{p+1}\setminus \{\Vect{0}\}$ and for all $\varphi_1,\ldots,\varphi_{p+1} \in A_1\times \cdots \times A_{p+1}$, the matrix $(p+1)\times (p+1)$ with entries $(\Mat{U})_{i,j}=u_j\Dp{0}{\varphi_i}$ (that does not depend on $\varphi_1,\ldots,\varphi_{p+1}$) is such that $\Mat{U} \Vect{y}\neq \Vect{0}$. So there exists $i_0(\Vect{y}) \in \{1,\ldots,p+1\}$ such that $\tr{\Vect{y}} \Mat{U}_{i_0(\Vect{y}),.} =\tr{\Vect{y}} \Vect{u}\Dp{0}{\varphi_{i_0(\Vect{y})}}\neq 0$. Therefore, for all $\Vect{y} \in \RR^{p+1}\setminus \{\Vect{0}\}$
$$P_{\ParVT} \Big( \Big\{ \varphi, \;\; \tr{\Vect{y}} \Vect{u}(0|\varphi)\neq 0 \Big\} \Big)> P_{\Vect{\ParVT}}(A_{i_0(\Vect{y})})>0,$$
which ends the proof.
\end{list}
\end{proof}

\section{Asymptotic normality for maximum pseudo-likelihood estimates} \label{sec-norm}

In this section, the existence of an ergodic measure is ensured, relatively to our framework, by Assumptions $\mathbf{E_1}$, $\mathbf{E_2^{loc}}$ and $\mathbf{E_3}$. 
The main tool used hereafter is a central limit theorem proposed by \cite{Jensen94} which justifies the need of $\mathbf{E_2^{loc}}$ instead of $\mathbf{E_2^{qloc}}$.

To ensure the asymptotic normality for the maximum pseudo-likelihood estimator, the following assumptions are needed. Denote, for some real $z$, by $[z]$ the integer part of $z$.
\begin{itemize}
\item[$\mathbf{N_1}$] The point process is observed in a domain $\Lambda_n\oplus D= \cup_{x \in \Lambda_n}\mathcal{B}(x,D)$, where $\Lambda_n\subset \RR^2$ can de decomposed into $\cup_{i \in I_n} \Dom{i}$ where for $i=(i_1,i_2)$
$$  
\Dom{i}=\left\{ z\in \RR^2, \Dt\left(i_j -\frac12\right) \leq z_j \leq \Dt\left(i_j -\frac12\right), j=1,2
\right\}$$ 
for some $\Dt>0$. As $n \to +\infty$, we also assume that $\Lambda_n \to \RR^2$ such that $|\Lambda_n|\to +\infty$ and ${\displaystyle \frac{|\partial \Lambda_n|}{|\Lambda_n|}\to 0}$

\item[$\mathbf{N_2}$] $\VIPar{0}{\cdot}{\ParV}$ is twice times differentiable in $\ParV=\ParVT$ and for all $j,k=1,\ldots,p+1$, there exists $\varepsilon>0$ such that the variables
$$
\dVIPar{j}{0}{\cdot}{\ParVT}^{3+\varepsilon} \; \mbox{ and } \; \ddVIPar{j}{k}{0}{\cdot}{\ParVT}\; \in L^1(P_{\ParVT})
$$
\item[$\mathbf{N_3}$] The matrix 
\begin{equation} \label{eq-defSig}
\Mat{\Sigma} ( \Dt,\ParVT ) = \Dt^{-2} \sum_{|i|\leq \left[ \frac{D}{\Dt}\right]+1} \Esp_{\ParVT} \left(
\dLPL{0}{\Phi}{\ParVT} \tr{\dLPL{i}{\Phi}{\ParVT} }
\right)
\end{equation}
is symmetric and definite positive. The vector $\dLPL{i}{\varphi}{\ParV}$ is defined  for any finite configuration $\varphi$ and for all $\Vect{\theta} \in \SpPar$ and $j=1,\ldots,p+1$ by
$$
\left( \dLPL{i}{\varphi}{\ParV} \right)_j = \int_{\Dom{i}} 
\dVIPar{j}{x}{\varphi}{\ParV}
\exp \left(-\VIPar{x}{\varphi}{\ParV}\right) dx - \sum_{x \in \varphi_{\Dom{i}} }
\dVIPar{j}{x}{\varphi \setminus x}{\ParV}.
$$

\item[$\mathbf{N_4}$] $\forall \Vect{y} \in \RR^{p+1}\setminus \{\Vect{0}\}$
$$
P_{\ParVT} \left( \Big\{  
\varphi, \; \tr{\Vect{y}} \Vect{V}^{(1)} (0|\varphi;\ParVT) \neq 0 
\Big\} \right) >0,
$$
where for $i=1,\ldots,p+1$, $(\Vect{V}^{(1)} (0|\varphi;\ParVT))_i=\dVIPar{i}{0}{\varphi}{\ParVT}$.

\item[$\mathbf{N_5}$] There exists a neighborhood $\Vois$ of $\ParVT$ such that $\VPar{\cdot}{\ParV}$ is twice times continuously differentiable for all $j,k=1,\ldots,p+1$, we have 
$$
\left|
\dVIPar{j}{0}{\Phi}{\ParV} - \dVIPar{j}{0}{\Phi}{\ParVT}
\right| \leq || \ParV -\ParVT ||^{c_1} \; h_1 (0,\Phi),
$$
and 
$$
\left|
\ddVIPar{j}{k}{0}{\Phi}{\ParV} - \ddVIPar{j}{k}{0}{\Phi}{\ParVT}
\right| \leq || \ParV -\ParVT ||^{c_2} \; h_2 (0,\Phi),
$$
with $c_1,c_2>0$ and $h_1(\cdot,\cdot),h_2(\cdot,\cdot)$ two functions such that, for all $x$, $h_i(0,\Phi)=h_i(x,\Phi_x)$ and such that $h_1(0,\cdot)^2$ and $h_2(0,\cdot) \in L^1(P_{\ParVT})$.
\end{itemize}

\begin{remark} Assumption $\mathbf{N_1}$ is similar to the one of \Citet{Jensen93}, \Citet{Jensen94} and \Citet{Heinrich92}. Among other things, $\mathbf{N_1}$ ensures that $\Lambda_n$ is a regular sequence of domains such that $\Lambda_n\to \RR^2$. 
\end{remark}

\begin{remark} \label{rem-NsansLS} Similarly to Remark~\ref{rem-CsansLS},
\begin{list}{$\bullet$}{}
\item the integrability condition occuring in Assumption $\mathbf{N_2}$ becomes~: $\big( \dVIPar{j}{0}{\cdot}{\ParVT} \exp \big( -\VIPar{0}{\cdot}{\ParVT} \big) \big)^{3+\varepsilon}$ and $\ddVIPar{j}{k}{0}{\cdot}{\ParVT} \exp \big( -\VIPar{0}{\cdot}{\ParVT} \big)$ are $P_{\ParVT}$-integrable.
\item the functions $h_1(\cdot,\cdot)$ and $h_2(\cdot,\cdot)$ occuring in Assumption $\mathbf{N_5}$ are now such that for all $\ParV \in \mathcal{V}$, the variables $h_1(0,\cdot)^2 \exp\big(- \VIPar{0}{\cdot}{\ParV}\big)$ and $h_2(0,\cdot)\exp\big( - \VIPar{0}{\cdot}{\ParV} \big)$ are $P_{\ParVT}$-integrable. Moreover, it is also assumed, for all $j,k=1,\ldots,p+1$, the $P_{\ParVT}$-integrability of the variables $\big( \dVIPar{j}{0}{\cdot}{\ParVT} \big)^2  \exp\big( -\VIPar{0}{\cdot}{\ParV} \big)$ and $\ddVIPar{j}{k}{0}{\cdot}{\ParVT}  \exp\big( -\VIPar{0}{\cdot}{\ParV} \big)$.
\end{list}
These Assumptions have been verified in \cite{Mase99} for the Ruelle class of pairwise interaction function with $\Vect{\theta}=(\beta,z)$ where $\beta$ represents the inverse temperature and $z$ the chemical potential.
\end{remark}

For $\ParV$ in a neighborhood of $\ParVT$, we can define, under Assumption $\mathbf{N_5}$, $\Vect{U}_n^{(1)}(\ParV)$ as the vector derivative of $U_n$. More precisely under Assumption $\mathbf{N_1}$, we can write
\begin{equation} \label{eq-Un1}
\Vect{U}_n^{(1)} (\ParV) = |\Lambda_n|^{-1}  \sum_{i \in I_n} \dLPL{i}{\varphi}{\ParV}.
\end{equation}
For $\ParV$ in a neighborhood of $\ParVT$, we can also define, under Assumption $\mathbf{N_5}$, the Hessian matrix $\Mat{U_n}^{(2)}(\ParV)$ given for $j,k=1,\ldots,p+1$ by
\begin{eqnarray}
\Mat{U_n}^{(2)}(\ParV) &=& \frac{1}{|\Lambda_n|} \int_{\Lambda_n} \left(
\dVIPar{j}{x}{\varphi}{\ParV} \dVIPar{k}{x}{\varphi}{\ParV}
- \ddVIPar{j}{k}{x}{\varphi}{\ParV}
\right) \exp\left( -\VIPar{x}{\varphi}{\ParV} \right)dx \nonumber \\
&& + \frac{1}{|\Lambda_n|}
\sum_{x \in \varphi_{\Lambda_n}} \ddVIPar{j}{k}{x}{\varphi\setminus x}{\ParV}.\label{eq-U2}
\end{eqnarray}

\begin{proposition} \label{prop-TCL}
Assume $P_{\ParVT}$ stationary, then under Assumptions $\mathbf{N_1}$ to $\mathbf{N_5}$, we have, for any $\Dt$ fixed, the following convergence in distribution as $n \to +\infty$
\begin{equation} \label{convLoiMPLE}
|\Lambda_n|^{1/2} \; \Estn{\Mat{\Sigma}}(\Dt,\Estn{\ParV})^{-1/2}
\; \Mat{U_n}^{(2)}(\Estn{\ParV})
 \; \left( \Estn{\ParV} - \ParVT  \right) \rightarrow 
\mathcal{N} \left( 0 , \Mat{I}_{p+1} \right),
\end{equation}
where for some $\ParV$ and some finite configuration $\varphi$, the matrix $\Estn{\Mat{\Sigma}}(\Dt,\ParV)$ is defined by
\begin{equation} \label{eq-defSigEst}
\Estn{\Mat{\Sigma}}(\Dt,\ParV)  = |\Lambda_n|^{-1} \Dt^{-2} \sum_{i \in I_n}\sum_{|j-i|\leq \left[ \frac{D}{\Dt}\right]+1,j\in I_n}   \dLPL{i}{\varphi}{\ParV} \tr{\dLPL{j}{\varphi}{\ParV} }
\end{equation}
\end{proposition}

By similar arguments of \Citet{Jensen94}, due to the decomposition of stationary measures as a mixture of ergodic measures (see \cite{Preston76}), one only needs to prove Proposition~\ref{prop-TCL} by assuming that $P_{\ParVT}$ is ergodic. Therefore, in Lemmas~\ref{lem-H2Guyon} to~\ref{lem-H3Guyon}, $P_{\ParVT}$ is assumed to be ergodic.

The proof of this result is based on a general result obtained by \Citet{Guyon92} (Proposition 3.7), giving conditions for which a miminum contrast estimator is asymptotically normal. The following Lemmas are needed to ensure these conditions. The first one ensures a central limit theorem for $\Vect{U}_n^{(1)}(\ParVT)$.

\begin{lemma} \label{lem-H2Guyon}
Under Assumptions $\mathbf{N_1}$, $\mathbf{N_2}$ and $\mathbf{N_3}$, \\
\noindent $(a)$ we have, for any fixed $\Dt$, the following convergence in distribution as $n \to +\infty$
\begin{equation} \label{Un1Loi}
|\Lambda_n|^{1/2} \Mat{\Sigma}(\Dt,\ParVT)^{-1/2}\; \Vect{U_n}^{(1)}(\ParVT) \rightarrow \mathcal{N} \left( 0, \Mat{I}_{p+1} \right)
\end{equation}
where the matrix $ \Mat{\Sigma}(\Dt,\ParVT)$ is defined by (\ref{eq-defSig}).

\noindent $(b)$ Moreover, we have $P_{\ParVT}-$almost surely as $n\to +\infty$
\begin{equation} \label{eq-convMatVn}
\Estn{\Mat{\Sigma}}(\Dt,\ParVT) \rightarrow \Mat{\Sigma}(\Dt,\ParVT).
\end{equation}
\end{lemma}

\begin{proof} $(a)$ The idea is to apply to $\Vect{U}_n^{(1)}(\ParVT)$ a central limit theorem obtained by \cite{Jensen94}, Theorem~2.1. The following conditions have to be fullfilled to apply this result.
\begin{itemize}
\item[$(i)$] For all $i\in I_n$ and for all $j=1,\ldots,p+1$ $\Esp_{P_{\ParVT}} \left( \left| (\dLPL{i}{\Phi}{\ParVT})_j \right|^3 \right) <+\infty.$
\item[$(ii)$] For all $i \in I_n$, $\Esp_{P_{\ParVT}} \left(
(\dLPL{i}{\Phi}{\ParVT} )_j | \Phi_{\Dom{i}^c} \right)=0.$
\item[$(iii)$] The set $I_n$ is such that $\dfrac{|\partial I_n|}{|I_n|} \to 0$, as $n \to +\infty.$
\item[$(iv)$] The matrix $\Var_{P_{\ParVT}} \left( |\Lambda_n|^{1/2} \Vect{U}_n^{(1)}(\ParVT) \right)$ converges to the matrix $\Mat{\Sigma}(\Dt,\ParV)$, which is definite positive under Assumption~$\mathbf{N_3}$.
\end{itemize}

\bigskip

\noindent \underline{Condition $(i)$}~: let us write
\begin{equation} \label{eq-T1T2}
\Esp_{P_{\ParVT}} \left( \left| (\dLPL{i}{\Phi}{\ParVT})_j
\right|^3 \right) \leq 2\times \left( T_1 +T_2 \right) 
\end{equation}
where the terms $T_1$ et $T_2$ are respectively defined by
\begin{eqnarray}
T_1 &=& \Esp_{P_{\ParVT}} \left( \left|
\int_{\Dom{i}} \dVIPar{j}{x}{\Phi}{\ParVT} \exp(-\VIPar{x}{\Phi}{\ParVT})dx
\right|^3 \right) \nonumber \\
T_2&=&  
\Esp_{P_{\ParVT}} \left( \left|
\sum_{x \in \Phi_{\Dom{i}}} \dVIPar{j}{x}{\Phi\setminus x}{\ParVT} \right|^3 \right) \nonumber 
\end{eqnarray}
Under Assumption $\mathbf{N_2}$, Hölder's inequality and the stationarity of $P_{\ParVT}$, we can prove that
\begin{eqnarray} 
T_1 &\leq& \exp(3K) |\Lambda_0|^2 \times \Esp_{P_{\ParVT}} \left( \int_{\Dom{i}} 
\left|  \dVIPar{j}{x}{\Phi}{\ParVT}\right|^3 dx \right) \nonumber \\
&\leq &\exp(3K) |\Lambda_0|^3 \times  \Esp_{P_{\ParVT}} \left(\left| \dVIPar{j}{0}{\Phi}{\ParVT}\right|^3 \right) <+\infty. \label{eq-T1}
\end{eqnarray}
We have from Hölder's inequality 
$$
T_2 \leq \Esp_{P_{\ParVT}} \left( |\Phi_{\Dom{i}}|^2 \sum_{x \in \Phi_{\Dom{i}}} \left| \dVIPar{j}{x}{\Phi\setminus x}{\ParVT}\right|^3 \right).
$$
And from Corollary~\ref{cor-CampGlotz}, it follows
$$T_2 \leq   \; |\Lambda_0| \;\Esp_{P_{\ParVT}}\left( |\Phi_{\Lambda_0}|^2 \left| \dVIPar{j}{0}{\Phi}{\ParVT}  \right|^3 \exp\left( -\VIPar{0}{\Phi}{\ParVT} \right)     \right).$$
Again from Hölder's inequality and under Assumption $\mathbf{E_3}$, we can prove that for all $\eta >0$,
$$
T_2 \leq \; |\Lambda_0| \exp(K)\; \Esp_{P_{\ParVT}} \Big( |\Phi_{\Lambda_0}|^{2(1+\frac1\eta)} \Big)^{\frac\eta{1+\eta}} \;\; 
\Esp_{P_{\ParVT}} \left( \left| \dVIPar{j}{0}{\Phi}{\ParVT}\right|^{3(1+\eta)} \right)^{\frac1{1+\eta}}.
$$
Under Assumption $\mathbf{E_3}$ it is well-known that, for all $z>0$, $\Esp_{P_{\ParVT}}(|\Phi_{\Lambda_0}|^z)<+\infty$. Now, let $\varepsilon=3\eta$, there exists $\kappa=\kappa(\varepsilon)$ such that
\begin{equation} \label{eq-T2}
T_2 \;\leq \; \kappa \; |\Lambda_0| \; \Esp_{P_{\ParVT}} \left( \left| \dVIPar{j}{0}{\Phi}{\ParVT}\right|^{3+\varepsilon} \right)^{\frac1{1+\varepsilon/3}} <+\infty
\end{equation}
under Assumption $\mathbf{N_2}$. Condition $(i)$ is obtained by combining~(\ref{eq-T1T2}), (\ref{eq-T1}) and~(\ref{eq-T2})

\noindent \underline{Condition $(ii)$}~: From the stationarity of the process, it is sufficient to prove that 
$$\Esp_{P_{\ParVT}}\left( (\dLPL{0}{\Phi}{\ParVT})_j | \Phi_{\Lambda_0^c} \right)=0.$$
Let us write for any finite configuration $\varphi$
\begin{equation} \label{goj}
 (\dLPL{0}{\varphi}{\ParVT})_j= -\int_{\Lambda_0} \dVIPar{j}{x}{\varphi}{\ParVT} 
 \exp(-\VIPar{x}{\varphi}{\ParVT})dx + 
\int_{\Lambda_0}\dVIPar{j}{x}{\varphi \setminus x}{\ParVT} \varphi(dx).
\end{equation}
Denote respectively by $G_1(\varphi)$ and $G_2(\varphi)$ the first and the second right-hand term of~(\ref{goj}) and by $E_i=\Esp_{P_{\ParVT}}\left( G_i(\Phi) | \Phi_{\Lambda_0^c}=\varphi_{\Lambda_0^c} \right)$. From the definition of Gibbs point processes,
$$
E_2 = \frac1{Z_{\Lambda_0}(\varphi_{\Lambda_0^c})} \int_{\Omega_{\Lambda_0}} Q(d\varphi_{\Lambda_0}) 
\int_{\RR^2} \varphi_{\Lambda_0}(dx)\mathbf{1}_{\Lambda_0}(x) \dVIPar{j}{x}{\varphi\setminus x}{\ParVT} \exp\left(-\VIPar{\varphi_{\Lambda_0}}{\varphi_{\Lambda_0^c}}{\ParVT}\right).
$$
Denote by $\varphi^\prime=(\varphi_{\Lambda_0},\varphi_{\Lambda_0^c}^\prime)$. Since $Q$ is a Poisson process we can write
{ \begin{eqnarray}
E_2&=& \frac1{Z_{\Lambda_0}(\varphi_{\Lambda_0^c})} \int_{\Omega} Q(d\varphi^\prime) 
\int_{\RR^2} \varphi^\prime(dx) \mathbf{1}_{\Lambda_0}(x) 
 \dVIPar{j}{x}{\varphi\setminus x}{\ParVT}
 \exp \left(-\VIPar{\varphi_{\Lambda_0}}{\varphi_{\Lambda_0^c}}{\ParVT} \right) \nonumber\\
&=&\frac1{Z_{\Lambda_0}(\varphi_{\Lambda_0^c})} \int_{\Omega} Q(d\varphi^\prime) 
\int_{\RR^2} \varphi^\prime(dx) \mathbf{1}_{\Lambda_0}(x) 
\dVIPar{j}{x}{\varphi_{\Lambda_0}^\prime \cup \varphi_{\Lambda_0^c} \setminus x}{\ParVT}
 \exp \left(-\VIPar{\varphi_{\Lambda_0}^\prime}{\varphi_{\Lambda_0^c}}{\ParVT} \right) \nonumber
\end{eqnarray} }
Now, from Campbell Theorem (applied to the Poisson measure $Q$)
$$E_2=\frac1{Z_{\Lambda_0}(\varphi_{\Lambda_0^c})} \int_{\Lambda_0}  dx \int_{\Omega} Q_x^!(d\varphi^\prime) 
\dVIPar{j}{x}{\varphi_{\Lambda_0}^\prime \cup \varphi_{\Lambda_0^c}}{\ParVT}
 \exp \left(-\VIPar{\varphi_{\Lambda_0}^\prime\cup x}{\varphi_{\Lambda_0^c}}{\ParVT} \right).
 $$
Since from Slivnyak-Mecke Theorem, $Q=Q_x^!$, one can obtain
\begin{eqnarray}
E_2 &=& \frac1{Z_{\Lambda_0}(\varphi_{\Lambda_0^c})} \int_{\Omega} Q(d\varphi^\prime) 
\int_{\Lambda_0}  dx \; 
\dVIPar{j}{x}{\varphi_{\Lambda_0}^\prime \cup \varphi_{\Lambda_0^c}}{\ParVT}
\exp \left(-\VIPar{\varphi_{\Lambda_0}^\prime\cup x}{\varphi_{\Lambda_0^c}}{\ParVT} \right)
 \nonumber \\
&=& \frac1{Z_{\Lambda_0}(\varphi_{\Lambda_0^c})} \int_{\Omega} Q(d\varphi_{\Lambda_0})
\int_{\Lambda_0} dx \dVIPar{j}{x}{\varphi}{\ParVT} 
\exp\left(- \VIPar{x}{\varphi}{\ParVT}\right) 
\exp\left(- \VIPar{\varphi_{\Lambda_0}}{\varphi_{\Lambda_0^c}}{\ParVT}\right) \nonumber \\
&=& -E_1 \nonumber 
\end{eqnarray}

\noindent \underline{Condition $(iii)$:} this condition is equivalent to Assumption $\mathbf{N_1}$.

\noindent \underline{Condition $(iv)$:} 
let us start by noting that the vector $\dLPL{i}{\varphi}{\ParVT}$ depends only on $\varphi_{\Dom{j}}$ for $j$ such that $|j-i|\leq \left[ \frac{D}{\Dt}\right]+1$. From~(\ref{eq-Un1}), we can obtain
\begin{eqnarray}
 \Var_{\ParVT} \left( |\Lambda_n|^{1/2} \Vect{U}_n^{(1)}(\ParVT)  \right)
&=& |\Lambda_n|^{-1} \Var_{P_{\ParVT}} \left( \dLPL{i}{\Phi}{\ParVT} \right) \nonumber \\
&=&|\Lambda_n|^{-1} \sum_{i,j \in I_n} \Esp_{P_{\ParVT}} 
\left( \dLPL{i}{\Phi}{\ParVT} \tr{ \dLPL{j}{\Phi}{\ParVT}} \right)  \nonumber \\
&=& |\Lambda_n|^{-1} \sum_{i \in I_n} \Bigg\{ 
\sum_{|j-i|\leq \left[ \frac{D}{\Dt}\right]+1, j\in I_n}   \Esp_{P_{\ParVT}} 
\left(  \dLPL{i}{\Phi}{\ParVT} \tr{ \dLPL{j}{\Phi}{\ParVT}}\right)  \nonumber \\
%\nonumber \\&&  
&& \quad  + \sum_{|j-i|> \left[ \frac{D}{\Dt} \right]+1, j\in I_n}   \Esp_{P_{\ParVT}} 
\left( \dLPL{i}{\Phi}{\ParVT} \tr{ \dLPL{j}{\Phi}{\ParVT}} \right) \Bigg\}. \nonumber
\end{eqnarray}
Let $j\in I_n$ such that $|j-i|>\left[\frac{D}{\Dt} \right]+1$, then using condition~$(ii)$
\begin{eqnarray}
\Esp_{P_{\ParVT}} 
\left(  \dLPL{i}{\Phi}{\ParVT} \tr{ \dLPL{j}{\Phi}{\ParVT}} \right) &=& \Esp_{P_{\ParVT}} \left( \Esp \left( \dLPL{i}{\Phi}{\ParVT} \tr{ \dLPL{j}{\Phi}{\ParVT}}| \dLPL{j}{\Phi}{\ParVT}  \right)\right) \nonumber \\
&=& \Esp_{P_{\ParVT}} 
\left(\Esp \left( \dLPL{i}{\Phi}{\ParVT} | \dLPL{j}{\Phi}{\ParVT} 
 \right) \tr{ \dLPL{j}{\Phi}{\ParVT}}  \right) \nonumber \\
&=& 0 \nonumber
\end{eqnarray}
Now, denote by $\widetilde{I}$ the following set
$$
\widetilde{I} = \left\{ 
k \in I_n, |k-i| \leq \left[ \frac{D}{\Dt} \right] +1, \forall i \in \partial I_n
\right\}
$$
and (for the sake of simplicity) by $E_{i,j}$ the following mean
$$
E_{i,j}= \Esp_{P_{\ParVT}} \left( \dLPL{i}{\Phi}{\ParVT}  \tr{\dLPL{j}{\Phi}{\ParVT}} \right). 
$$
From the stationarity of the process, we can write
\begin{eqnarray}
\Var_{P_{\ParVT}} \left( |\Lambda_n|^{1/2} \Vect{U}_n^{(1)}(\ParVT)  \right)
&=& |\Lambda_n|^{-1} \left( \sum_{i \in I_n\setminus \widetilde{I}} \;\;\sum_{|j-i|\leq \left[ \frac{D}{\Dt} \right]+1, j\in I_n} E_{i,j} \; + \; \sum_{i \in \widetilde{I}} \;\;\sum_{|j-i|\leq \left[ \frac{D}{\Dt} \right]+1, j\in I_n} E_{i,j} \right). \nonumber \\
&=& |I_n\setminus\widetilde{I}| \; |\Lambda_n|^{-1} \sum_{|i|\leq \left[ \frac{D}{\Dt} \right]+1} E_{0,i} \quad + \quad |\widetilde{I}| \; |\Lambda_n|^{-1} \sum_{|j-i_0| \leq \left[ \frac{D}{\Dt} \right]+1} E_{i_0,j}, \nonumber
\end{eqnarray}
for some $i_0 \in \partial I_n$. From the definition of the set $I_n$, we have as $n \to +\infty$
$$
\Var_{P_{\ParVT}} \left( |\Lambda_n|^{1/2} \Vect{U}_n^{(1)}(\ParVT)  \right) \to \sum_{|i| \leq \left[ \frac{D}{\Dt} \right]+1} E_{0,i} \; = \; \Mat{\Sigma}\left( \Dt , \ParVT \right).
$$

$(b)$ According to~(\ref{eq-defSigEst}), it is easy to see that $\Estn{\Mat{\Sigma}}(\Dt,\ParVT)$ is defined such that as $n \to +\infty$, $$\Esp_{P_{\ParVT}}\left( \Estn{\Mat{\Sigma}}(\Dt,\ParVT) \right) \to \Mat{\Sigma}(\Dt,\ParVT).$$
We leave the reader to check that under Assumption $\mathbf{N_1}$ and from Theorem~1 (\cite{Nguyen79}), we have $P_{\ParVT}-$almost surely as $n \to +\infty$, $\Estn{\Mat{\Sigma}}(\Dt,\ParVT) \to \Mat{\Sigma}(\Dt,\ParVT)$.
\end{proof}

\begin{remark} \label{rem-hypN3}
From the previous proof, we can note that Assumption $\mathbf{N_3}$ is fullfilled as soon as one can prove that for $n$ sufficiently large the matrix $\Var_{P_{\ParVT}} \left( |\Lambda_n|^{1/2} \Vect{U}_n^{(1)}(\ParVT)  \right)$ is definite positive.
\end{remark}

\begin{lemma} \label{lem-H1Guyon}
Under Assumptions $\mathbf{N_1}$, $\mathbf{N_2}$ and $\mathbf{N_5}$, there exists a neighborhood $\Vois$ of $\ParVT$ on which $U_n(\cdot)$ is twice times continuously differentiable and a random variable $Y$ such that for all $j,k=1,\ldots,p+1$ and for all $\ParV \in \Vois$ we have,  
$$\Big| \left(\Mat{U_n}^{(2)}(\ParV) \right)_{j,k}\Big| \leq Y.$$
\end{lemma}

\begin{proof}
Let $j,k=1,\ldots,p+1$. Under Assmuption $\mathbf{N_5}$, there exists a neighborhood $\Vois$ of $\ParVT$ such that we can write for any configuration $\varphi$
\begin{eqnarray}
\left(\Mat{U_n}^{(2)}(\ParV) \right)_{j,k} &=& -\frac1{|\Lambda_n|} \int_{\Lambda_n} 
\ddVIPar{j}{k}{x}{\varphi}{\ParV} 
\exp\left( - \VIPar{x}{\varphi}{\ParV}\right)dx \nonumber \\ 
&& +  
\frac1{|\Lambda_n|} \int_{\Lambda_n} 
\dVIPar{j}{x}{\varphi}{\ParV} \dVIPar{k}{x}{\varphi}{\ParV}
\exp\left( - \VIPar{x}{\varphi}{\ParV}\right)dx \nonumber \\
&& + \frac1{|\Lambda_n|} \sum_{x \in \varphi_{\Lambda_n}} 
\ddVIPar{j}{k}{x}{\varphi\setminus x}{\ParV} .
\end{eqnarray}
Denote respectively by $R_1,R_2,R_3$ the three right-hand terms of the previous equation. Under Assumption $\mathbf{N_5}$, one can choose the neighborhood $\Vois$ such that $|| \ParV -\ParVT||\leq \kappa$. Thus, one can obtain
\begin{eqnarray}
|R_1| & \leq & \exp(K) \frac1{|\Lambda_n|} \int_{\Lambda_n} \left( \kappa^{c_2} h_2(x,\varphi) + \left| \ddVIPar{j}{k}{x}{\varphi}{\ParVT} \right| \right)dx \nonumber \\
|R_2| & \leq & \exp(K) \frac1{|\Lambda_n|} \int_{\Lambda_n} \Big(
\kappa^{2c_1} h_1(x,\varphi)^2 
+ \kappa^{c_1} h_1(x,\varphi) \left|\dVIPar{j}{x}{\varphi}{\ParVT}\right| \nonumber \\
&& +\kappa^{c_1} h_1(x,\varphi) \left|\dVIPar{k}{x}{\varphi}{\ParVT}\right|
+ \Big| \dVIPar{j}{x}{\varphi}{\ParVT}  \dVIPar{k}{x}{\varphi}{\ParVT}  \Big| 
\Big) \nonumber \\
|R_3| & \leq & \frac1{|\Lambda_n|} \sum_{x \in \varphi_{\Lambda_n}} \left( \kappa^{c_2} h_2(x,\varphi \setminus x) + \Big| \ddVIPar{j}{k}{x}{\varphi\setminus x}{\ParVT}  \Big| \right) \nonumber
\end{eqnarray}
Under Assumptions $\mathbf{N_1}$ and $\mathbf{N_2}$, from Theorem~\ref{thmNguyen} (\cite{Nguyen79}), and using the stationarity of $P_{\ParVT}$, there exists $N_0\in \NN$ such that for all $n \geq N_0$, we have $P_{\ParVT}-$almost surely
\begin{eqnarray}
|R_1| &\leq & 2\times \exp(K) \Esp_{P_{\ParVT}} \left( \kappa^{c_2} h_2(0,\Phi) + \left| \ddVIPar{j}{k}{0}{\Phi}{\ParVT} \right|\right) \nonumber \\
|R_2| &\leq & 2\times \exp(K) \left\{ \Esp_{P_{\ParVT}} \left( k^{2c_1} h_1(0,\Phi)^2 +
\Big| \dVIPar{j}{0}{\Phi}{\ParVT}  \dVIPar{k}{0}{\Phi}{\ParVT}  \Big| \right) \right. \nonumber \\
&& \left. + \Esp_{P_{\ParVT}} \left( \kappa^{c_1} h_1(0,\Phi) \Big|\dVIPar{j}{0}{\Phi}{\ParVT}\Big|+
\kappa^{c_1} h_1(0,\Phi) \Big|\dVIPar{k}{0}{\Phi}{\ParVT}\Big|
 \right) \right\} \nonumber \\
|R_3| &\leq & 2\times \exp(K) \Esp_{P_{\ParVT}} \left( \kappa^{c_2} h_2(0,\Phi) + \Big| \ddVIPar{j}{k}{0}{\Phi}{\ParVT}  \Big|\right) \nonumber
\end{eqnarray}
Consequenlty, for $n$ large enough, there exists a positive constant $\kappa^\prime$ such that
$\Big|\Big(\Mat{U_n}^{(2)}(\ParV)\Big)_{j,k}\Big| \leq \kappa^\prime$, which implies the result.
\end{proof}

\begin{lemma} \label{lem-H3Guyon}
Under Assumptions $\mathbf{N_1}$ and $\mathbf{N_2}$, we have almost surely, as $n \to +\infty$
$$
\Mat{U_n}^{(2)}(\ParVT) \rightarrow \Mat{U}^{(2)}(\ParVT)$$
where $\Mat{U}^{(2)}(\ParVT)$ is the $(p+1)\times (p+1)$ matrix whose entry is
\begin{equation} \label{eq-U2}
 \left(\Mat{U}^{(2)}(\ParVT) \right)_{j,k}  =
\Esp_{P_{\ParVT}} \left( 
\dVIPar{j}{0}{\Phi}{\ParVT} \dVIPar{k}{0}{\Phi}{\ParVT}
\exp\left( -\VIPar{0}{\Phi}{\ParVT} \right)\right).
\end{equation}
Furthermore, under Assumption $\mathbf{N_4}$, $\Mat{U}^{(2)}$ is a symmetric definite positive matrix.
\end{lemma}

\begin{proof}
Let $j,k=1,\ldots,p+1$. Under Assumptions $\mathbf{N_1}$ and $\mathbf{N_2}$ and from Theorem~\ref{thmNguyen} (\cite{Nguyen79}), we have almost surely, as $n \to +\infty$
\begin{eqnarray}
\left(\Mat{U_n}^{(2)}(\ParVT) \right)_{j,k} &\to & -\frac1{|\Lambda_0|} \Esp_{P_{\ParVT}} \left(  \int_{\Lambda_0} \ddVIPar{j}{k}{x}{\Phi}{\ParVT}
\exp\left( - \VIPar{x}{\Phi}{\ParVT}\right)dx \right)\nonumber \\ 
&& +  
\frac1{|\Lambda_0|} \Esp_{P_{\ParVT}} \left( 
\int_{\Lambda_0} 
\dVIPar{j}{x}{\Phi}{\ParVT} \dVIPar{k}{x}{\Phi}{\ParVT}
\exp\left( - \VIPar{x}{\Phi}{\ParVT}\right)dx \right)\nonumber \\
&& + \frac1{|\Lambda_0|} \Esp_{P_{\ParVT}} \left( 
\sum_{x \in \varphi_{\Lambda_0}} \ddVIPar{j}{k}{x}{\Phi \setminus x}{\ParVT}
\right)
\end{eqnarray}
Equation~(\ref{eq-U2}) is obtained using Corollary~\ref{cor-CampGlotz}. And under Assumption $\mathbf{N_4}$, it is easy to see that $\Mat{U}^{(2)}$ is a symmetric definite positive matrix.
\end{proof}

\bigskip \bigskip \bigskip \bigskip 

\noindent {\bf Proof of Proposition~\ref{prop-TCL}}
Using Lemmas~\ref{lem-H2Guyon} à~\ref{lem-H3Guyon}, one can apply a classical result concerning asymptotic normality for minimum contrast estimators, {\it e.g.} Proposition 3.7 de \Citet{Guyon92}, in order to prove as $n \to +\infty$ 
$$
|\Lambda_n|^{1/2} \Estn{\Mat{\Sigma}}(\Dt,\ParVT)^{-1/2} \; \Mat{U}_n^{(2)}(\ParVT) 
\; \left( \Estn{\ParV} - \ParVT  \right) \rightarrow 
\mathcal{N} \left( 0 , \Mat{I}_{p+1} \right).$$
The result is then obtained using the fact that $\Estn{\ParV}$ is a consistent estimator of $\ParVT$.
$\blacksquare$

Let us precise, as in Section~\ref{sec-cons}, the different Assumptions for energy functions that can be written as~(\ref{eq-modExp}).

\begin{proposition}
For energy functions described by~(\ref{eq-modExp}), Assumptions $\mathbf{N_2}$ and $\mathbf{N_5}$ (resp. $\mathbf{N_4}$) can be replaced by $\mathbf{N_{2,5}^{exp}}$ (resp. $\mathbf{N_4^{exp}}$)
\begin{itemize}
\item[$\mathbf{N_{2,5}^{exp}}$] For $i=1,\ldots,p+1, \quad $ there exists $\varepsilon>0$ such that ${\displaystyle u_i(0|\cdot) \in L^{3+\varepsilon}(P_{\ParVT}).}$
\item[$\mathbf{N_{4}^{exp}}= \mathbf{C_3^{exp}}$]
\end{itemize}
\end{proposition}

The proof is trivial.

\section{Some examples} \label{sec-ex}

In this section, it is assumed that the sequence of domains satisfies $\mathbf{N_1}$ (which implies $\mathbf{C_1}$). Moreover, we only focus on examples satisfying the following convenient Assumption denoted by $\mathbf{M}$~:
\begin{itemize}
\item[$\mathbf{M}$] \quad There exists $K_1,K_2>0$ such that for any finite configuration $\varphi$, we have for all $x$
$$
-K_1 \leq  u_i(x|\varphi)  \leq K_2,  \qquad \mbox{ for } i=1,\ldots,p+1.
$$
\end{itemize}
Quite obviously, Assumption $\mathbf{M}$ ensures $\mathbf{C_{2,4}^{exp}}$ and $\mathbf{N_{2,5}^{exp}}$. Let us now present a Corollary of Propositions~\ref{prop-convMPLE} and~\ref{prop-TCL}.

\begin{corollary}
Under Assumption $\mathbf{M}$ and $\mathbf{C_3^{exp}}$, the consistency of the maximum pseudo-likelihood, that is the result~(\ref{eq-convMPLE}), is valid. And in addition with $\mathbf{N_3^{exp}}$, its asymptotic normality property, that is the result~(\ref{convLoiMPLE}), is ensured.
\end{corollary}

\subsection{Pairwise $\beta$-Delaunay model}

We first deal with our main example: $\beta$-Delaunay (of order some small enough fixed $\beta_0$) model with multi-Strauss pairwise interaction function. In other words,
\begin{equation} \label{eq-MSDel2}
\VPar{\varphi}{\ParV}=\theta^{(1)}|\varphi|+\sum_{\xi\in Del_{2,\beta}^{\beta_0}(\varphi)} u^{(2)}(\xi;\varphi,\ParVi{2})=\tr\ParV \Vect{u}(\varphi)
\end{equation}
with $u_1(\varphi)=|\varphi|$ and for any $i\in\{2,\ldots,p+1\}$,
\[
u_i(\varphi)=\sum_{\xi\in Del_{2,\beta}^{\beta_0}(\varphi)} \mathbf{1}_{]d_{i-1},d_i]}(\|\xi\|)
\] 
where $0=d_1\leq d_2\leq \ldots\leq d_{p+1}$ are some fixed real numbers. Literally, $u_i(\varphi)$ ($i>1$) corresponds to the number of ($\beta$-Delaunay) edges of length between $d_{i-1}$ and $d_i$.
We may also notice that the range of the pairwise interaction function is $d_{p+1}$, that is $u^{(2)}(\xi;\varphi,\ParVi{2})=0$ when $\|\xi\|>d_{p+1}$.

In \cite{BBD4}, it is proved that this model satisfies Assumption $\mathbf{M}$. Let us now verify the technical conditions $\mathbf{C_3^{exp}}$ and $\mathbf{N_3^{exp}}$.

\begin{proposition} \label{prop-C3exp}
Assumption $\mathbf{C_3^{exp}}$ is satisfied for the $\beta$-Delaunay model with multi-Strauss pairwise interaction function.
\end{proposition}

\begin{proof}
Denote by $\Delta$ the following domain
$$
\Delta = \left\{ z \in \RR^2: -D\leq z_i\leq D, i=1,2\right\}
$$
and by $A_1$ the event $A_1=\{ \varphi, \varphi_\Delta=\emptyset\}$. We clearly have for all $\varphi_1 \in A_1$, $\Vect{u}(0|\varphi_1)=\tr{(1,0,\ldots,0)}$. Now, let us give for $j=2,\ldots,p+1$, the points $c_{1,j}$ and $c_{2,j}$ such that the distances $d(0,c_{1,j})=d(0,c_{2,j})=d(c_{1,j},c_{2,j})=\frac{d_{j-1}+d_j}{2}$. Denote for $j=2,\ldots,p+1$ the following events for some $\eta>0$
$$
A_j(\eta)= \Big\{  \varphi \in \Omega: \varphi_\Delta=\{z_{1},z_{2}\}, z_{1}\in \mathcal{B}(c_{1,j},\eta), z_{2,j}\in \mathcal{B}(c_{2,j},\eta)  \Big\}.
$$ 
One can choose $\eta$ such that for all $\varphi \in A_j(\eta)$, the distances $d(0,z_{1}), d(0,z_{2}$ and $d(z_{1},z_{2})$ are comprised between $d_{j-1}$ and $d_j$. One can also choose $\eta$ such that the smallest angle of the triangle with vertices $\{0,c_{1,j},c_{2,j}\}$ is strictly greater thant $\beta_0$, which means that $\{0,c_{1,j},c_{2,j}\} \in Del_3^{\beta_0}(\varphi_j)$. Now, it is easy to see that the matrix $\Mat{U}$ defined in Proposition~\ref{prop-C234exp} is given by
$$
\Mat{U}=(u_j(0|\varphi_i))_{1\leq i,j\leq p+1} = \left( 
\begin{array}{ccccc}
1& 0& \cdots &\cdots & 0 \\
1& 3 & \ddots & \vdots & \vdots \\
1 & 0 & \ddots & \ddots& \vdots \\ 
\vdots & \vdots & \ddots & 3 & 0 \\
1 & 0 & \ldots & 0 & 3 \\
\end{array} \right)
$$
and is clearly invertible, which ends the proof.
\end{proof}

\begin{proposition} \label{prop-N3} Assumption $\mathbf{N_3}$ is satisfied for the $\beta$-Delaunay model with multi-Strauss pairwise interaction function.
\end{proposition}

\begin{proof}
From Remark~\ref{rem-hypN3}, it is sufficient to prove that the matrix $\Var_{P_{\ParVT}}(|\Lambda_n|^{1/2} \Vect{U}_n^{(1)}(\Phi;\ParVT))$ is definite positive for $n$ sufficiently large. Let $\Dt>D$, $\Vect{y}\in \RR^{p+1}$ and let $\widetilde{\Lambda}=\cup_{|i|\leq1}\Dom{i}$, by the same argument of \cite{Jensen94} (Equation (3.2)), we can write
$$\tr{\Vect{y}} \Var_{P_{\ParVT}}\left( |\Lambda_n|^{1/2} \Vect{U}_n^{(1)}\right) \Vect{y} \geq
|\Lambda_n| \; \Esp_{P_{\ParVT}} \left( \Var_{P_{\ParVT}} \left(
\tr{\Vect{y}} \Vect{U}_n^{(1)} | \Phi_{\Lambda_\ell}, \ell \notin 3\ZZ^2
\right)\right).$$
Now, following the proof of Lemma~\ref{lem-H2Guyon} ($(a)$ condition $(iv)$), one can prove that there exists $n_0\in \NN$ such that for all $n\geq n_0$,
$$\tr{\Vect{y}} \Var_{P_{\ParVT}}\left( |\Lambda_n|^{1/2}\right) \Vect{y} \geq \frac12 \Esp_{P_{\ParVT}} \left(
\Var_{P_{\ParVT}}\left( 
\tr{\Vect{y}} \dLPLt{\Phi}{\ParVT} | \Phi_{\Lambda_{\ell}}, 1\leq |\ell|\leq 2
\right) \right).
$$
The aim is to prove that the function $h(\Phi)= \tr{\Vect{y}}\dLPLt{\Phi}{\ParVT}$ is not almost surely a constant, when the variables $\Phi_{\Lambda_{\ell}}, 1\leq |\ell|\leq 2$ are (for example) fixed to $\emptyset$. Assume that the function $h(\cdot)$ explicitly given for any finite configuration $\varphi$ by 
$$
h(\varphi)= \sum_{k=1}^{p+1} y_k \left\{
\int_{\widetilde{\Lambda}} u_k(x|\varphi_{\Lambda_0})\exp\left( -\tr{\ParVT}\Vect{u}(x|\varphi_{\Lambda_0}) \right)dx \;\; - \;\;
\sum_{x \in \varphi_{\Lambda_0}} u_k(x|\varphi_{\Lambda_0}\setminus x)
\right\}
$$
is constant for all $\varphi \in \Omega^\prime=\{ \varphi\in \Omega: \varphi_{\Lambda_{\ell}}=\emptyset, 1\leq |\ell|\leq 2\}$.

Denote by $A_0=\left\{ \varphi \in \Omega^\prime : \varphi_{\Lambda_0}=\emptyset\right\}$ and by $A_1=\left\{ \varphi \in \Omega^\prime : |\varphi_{\Lambda_0}|=1 \right\}$. It is clear that, $P_{\ParVT}(A_0)>0$ and $P_{\ParVT}(A_1)>0$. We have for all $\varphi_0\in A_0$ and for all $\varphi_1 \in A_1$
$$
h(\varphi_0) = y_1 |\widetilde{\Lambda}| \exp(-\theta_1^\star) \quad \mbox{ and } \quad
h(\varphi_1) = y_1 |\widetilde{\Lambda}| \exp(-\theta_1^\star)-y_1.
$$
Assuming $h(\cdot)$ constant implies that $y_1=0$ and then $h(\cdot)$ vanishes. We now consider particular configurations of two points in $\Lambda_0$ and empty in $\Lambda_\ell\setminus\Lambda_0, 1\leq |\ell|\leq 2$.
Let us first introduce the following sets for any $j\in \{1,\cdots,p-1\}$ and any $\eta>0$
\begin{eqnarray}
D_j(\eta)&=&\Big\{(z_1,z_2)\in \Lambda_0^2 : z_1\in\mathcal{B}\left((0,0),\frac{\eta}4\right) \mbox{ and }z_2\in \mathcal{B}\left((d_j,0),\frac{3\eta}4\right)\Big\} \nonumber\\
D_j^-(\eta)&=&\Big\{(z_1,z_2)\in \Lambda_0^2 : z_1\in\mathcal{B}\left((0,0),\frac{\eta}4\right)\mbox{ and }z_2\in \mathcal{B}\left((d_j-\frac{\eta}2,0),\frac{\eta}4\right)\Big\}\subset D_j(\eta)\nonumber\\
D_j^+(\eta)&=&\Big\{(z_1,z_2)\in \Lambda_0^2 : z_1\in\mathcal{B}\left((0,0),\frac{\eta}4\right)\mbox{ and }z_2\in \mathcal{B}\left((d_j+\frac{\eta}2,0),\frac{\eta}4\right)\Big\}\subset D_j(\eta) \nonumber
\end{eqnarray}
When $\eta$ is small enough, the couple of points $(z_1, z_2)\in D_j(\eta)$ (resp. $D_j^-(\eta)$ and $D_j^+(\eta)$) are such that
$d_{j-1}<d_j-\eta<d(z_1,z_2)<d_j+\eta<d_{j+1}$ (resp. $d_{j-1}<d_j-\eta<d(z_1,z_2)<d_j$ and  $d_j<d(z_1,z_2)<d_j+\eta<d_{j+1}$).

We now derive the corresponding events for any $j\in \{1,\cdots,p-1\}$ and any $\eta>0$
\begin{eqnarray}
A_j(\eta) &=& \Big\{ \varphi \in \Omega^\prime: \varphi_{\Lambda_0}=\{z_1,z_2\} \mbox{ with }(z_1,z_2)\in D_j(\eta) \Big\} \nonumber \\
A_j^-(\eta) &=& \Big\{ \varphi \in \Omega^\prime : \varphi_{\Lambda_0}=\{z_1,z_2\} \mbox{ with }(z_1,z_2)\in D_j^-(\eta) \Big\} \subset A_j(\eta)\nonumber \\
A_j^+(\eta) &=& \Big\{ \varphi \in \Omega^\prime : \varphi_{\Lambda_0}=\{z_1,z_2\} \mbox{ with }(z_1,z_2)\in D_j^+(\eta) \Big\} \subset A_j(\eta) \nonumber 
\end{eqnarray}
satisfying $P_{\ParVT}(A_j(\eta) )>0$,  $P_{\ParVT}(A_j^-(\eta) )>0$ and $P_{\ParVT}(A_j^+(\eta) )>0$. 

Let us fix some $\varphi\in A_j(\eta)$. There exists some unique couple of points $(z_1,z_2)\in D_j(\eta)$, for which we define the following domain
\[
\widetilde{\Lambda}(z_1,z_2)= \left\{  x \in \widetilde{\Lambda}: \{x,z_1,z_2\} \in Del_3^{\beta_0}(\varphi\cup\{x\}) \right\}.
\]
Since $\{z_1,z_2\}\notin Del_2^{\beta_0}(\varphi)$, we then derive that
\begin{equation}\label{eq-hphi2}
0=h(\varphi)=\sum_{k=2}^{p+1} y_k 
\int_{\widetilde{\Lambda}(z_1,z_2)} u_k(x|\{z_1,z_2\})\exp\left( -\tr{\ParVT}\Vect{u}(x|\{z_1,z_2\}) \right)dx.
\end{equation}
When $\{x,z_1,z_2\} \in Del_3^{\beta_0}(\varphi\cup\{x\})$, we decompose $u_k(x|\{z_1,z_2\})$ into two additive terms in order to isolate the contribution of $x$: 
\[
u_k(x|\{z_1,z_2\})=u_k^{-x}(z_1,z_2)+u_k^x(z_1,z_2)
\]
with $u_1^{-x}(z_1,z_2)=0$ and $u_1^{x}(z_1,z_2)=1$, and for $k\neq 1$, 
$u_k^{-x}(z_1,z_2)=\mathbf{1}_{[d_{k-1},d_{k}[}(\|z_1-z_2\|)$ and  $u_k^{x}(z_1,z_2)=\displaystyle{\sum_{j=1,2}\mathbf{1}_{[d_{k-1},d_{k}[}(\|x-z_j\|)}$.\\
Then Equation~(\ref{eq-hphi2}) becomes for any $\Vect{y}\in\RR^{p+1}$ with $y_1=0$
\begin{eqnarray}
&&h(\varphi)=0= \exp\left(-\tr{\ParVT}\Vect{u}^{-x}(z_1,z_2) \right)\left(\tr{\Vect{y}} \Vect{u}^{-x}(z_1,z_2) f_1(z_1,z_2)+\tr{\Vect{y}} \Vect{f}(z_1,z_2)\right)\nonumber\\
&\Longleftrightarrow& \tr{\Vect{y}} \Vect{u}^{-x}(z_1,z_2) f_1(z_1,z_2)+\tr{\Vect{y}} \Vect{f}(z_1,z_2)=0 \label{eq-hphi2bis}
\end{eqnarray}
where $\Vect{u}^{-x}=(u_1^{-x},\cdots,u_{p+1}^{-x})$, $\Vect{u}^x=(u_1^x,\cdots,u_{p+1}^x)$ and $\Vect{f}=(f_1,f_2,\ldots,f_{p+1})$ with
\[
f_k(z_1,z_2) = \int_{\widetilde{\Lambda}(z_1,z_2)}  u_k^x(z_1,z_2)\exp\left(-\tr{\ParVT}\Vect{u}^x(z_1,z_2) \right) dx 
\]
Since each $f_k$ is continuous, one could assert that for any $\varepsilon>0$, there exists $\eta>0$ such that for any $(z_1,z_2)\in D_j(\eta)$, $|f_k(z_1,z_2)-\widetilde{f}_k|<\varepsilon$ where $\widetilde{f}_k=f_k((0,0),(d_j,0))$ is positive. We then set $\delta_k(z_1,z_2)=f_k(z_1,z_2)-\widetilde{f}_k$ and $\Vect{\delta}=(\delta_1,\cdots,\delta_{p+1})$.
We now apply the equation~(\ref{eq-hphi2bis}) for some fixed $\varphi_j^-\in A_j^-(\eta)$ and $\varphi_j^+\in A_j^+(\eta)$.  By denoting $(z_1^-,z_2^-)\in D_j^-(\eta)$ and $(z_1^+,z_2^+)\in D_j^+(\eta)$ such that $\varphi_j^-\cap\Lambda_0=\{z_1^-,z_2^-\}$ and $\varphi_j^+\cap\Lambda_0=\{z_1^+,z_2^+\}$, we have
\begin{eqnarray}
y_{j} f_1(z_1^-,z_2^-)+\tr{\Vect{y}} \Vect{f}(z_1^-,z_2^-)=0\nonumber\\
y_{j+1} f_1(z_1^+,z_2^+)+\tr{\Vect{y}} \Vect{f}(z_1^+,z_2^+)=0\nonumber
\end{eqnarray}
By substracting these two terms, we can obtain
\begin{equation} \label{eq-condFinale}
(y_{j+1}-y_{j})\widetilde{f}_1 = y_{j} \delta_1(z_1^-,z_2^-)-
y_{j+1} \delta_1(z_1^+,z_2^+) + \tr{\Vect{y}} \left(
{\Vect{\delta}}(z_1^-,z_2^-) 
- \Vect{\delta}(z_1^+,z_2^+) \right)
\end{equation}
By the previous continuity argument on the $f_k$, on can choose $\eta>0$ (depending on $\Vect{y}$) small enough such that the absolute value of the right-hand term of (\ref{eq-condFinale}) could be lower than any $\varepsilon>0$. Thus, by assuming that $y_j\neq y_{j+1}$ and choosing $\varepsilon=\frac12|y_{j+1}-y_{j}|\widetilde{f}_1$, there exists $\eta$ such that $|y_{j+1}-y_{j}| \widetilde{f}_1\leq \frac12 |y_{j+1}-y_{j}| \widetilde{f}_1$ which leads to an obvious contradiction.  Thus, (\ref{eq-condFinale}) holds only if $y_{j+1}=y_{j}$. By iterating this argument, we obtain that $y_2=y_3=\ldots=y_{p+1}$ and by applying this result on the equation~(\ref{eq-hphi2bis}), one may assert that for any $y\in\RR$ 
\[
y \sum_{k=1}^{p+1}f_k(z_1,z_2)= 0 ,
\] 
which implies that $y=0$ since $f_k(z_1,z_2)$ is positive for any $(z_1,z_2)\in D_j(\eta)$ ($j$ arbitrarily chosen in $\{1,\cdots,p+1\}$)
\end{proof}

\begin{figure}
\begin{tabular}{cc}
\includegraphics[width=8cm,height=9cm]{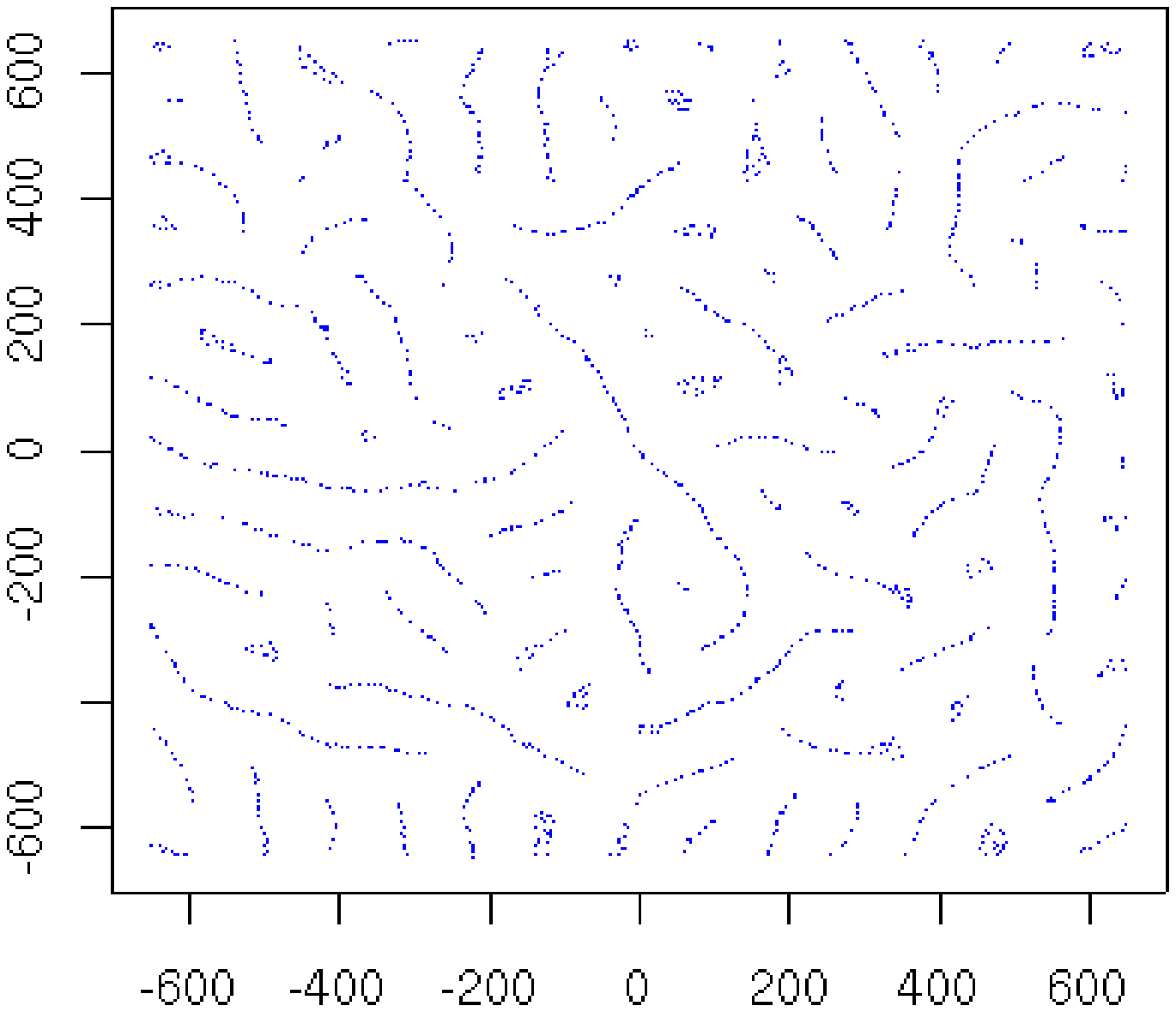} &
\includegraphics[width=8cm,height=9cm]{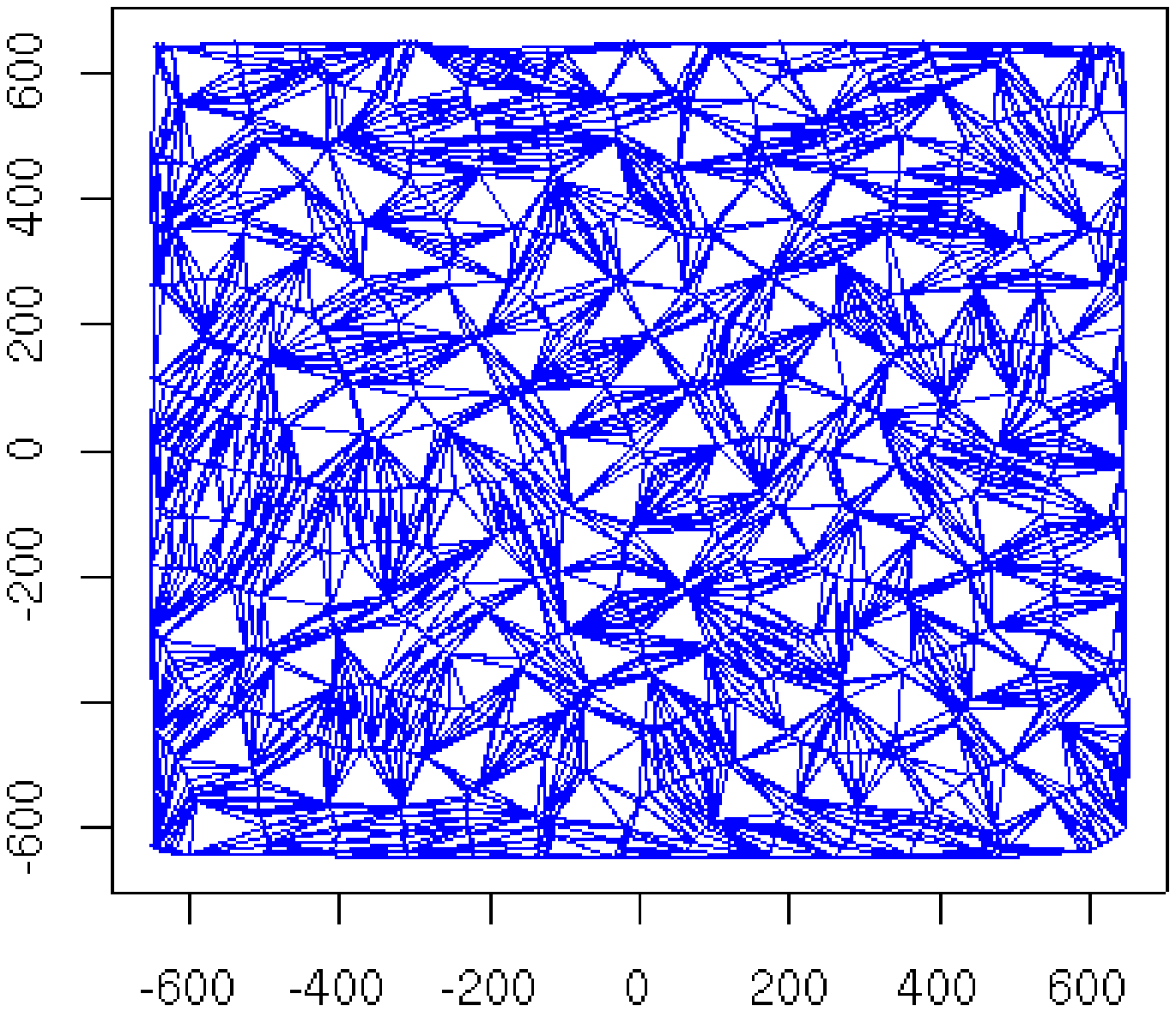} \\
\end{tabular}
\caption{Example of the points, with (left) and without the edges (right), of the realization of the $\beta$-Delaunay model with multi-Strauss pairwise interaction function where parameters $\ParV$ and $\Vect{d}$ are respectively fixed to $(0,2,4)$ and $(0,20,80)$.}
\label{fig-exDelPP}
\end{figure}

We propose a simulation study to verify the consistency of maximum pseudo-likelihood estimator. We consider the model~(\ref{eq-MSDel2}) with the vector of parameters $\Vect{\theta}=(0,2,4)$. The vector of bounds $\Vect{d}$ is assumed to be known and fixed to $\Vect{d}=(0,20,80)$. The simulation procedure used here is a direct adaptation to the Delaunay energies of the Geyer and M{{\o}}ller proposal~(\cite{Geyer94}, \cite{Geyer99}). We refer the reader to~\cite{BBD1} for a detail of the used algorithm. One simulation of such a point process is proposed in Figure~\ref{fig-exDelPP}. Table~\ref{tab-res} summarizes the different results obtained via $m=5000$ replications each one is generated after one million of iterations of the algorithm. One may verify that both the bias and the standard deviation become smaller and smaller as the domain $\Lambda_n$ grows.

\subsection{Other examples of pairwise interaction models}

In order to satisfy $\mathbf{C_{3}^{exp}}$ and $\mathbf{N_{3}^{exp}}$ for models on the complete graph or on the $k$ nearest-neighbours graph with multi-Strauss pairwise interaction function, we can chosse as in~\cite{Jensen94} a configuration with one point or two points. On the delaunay graph, it may be interesting to study multi-Strauss interaction function on the circumradius or on the smallest angle of each Delaunay triangle. As discussed previously the identifiability assumption $\mathbf{C_{3}^{exp}}$ holds easily but $\mathbf{N_{3}^{exp}}$ needs more attention. Otherwise, for pairwise Delaunay model, we can replace the assumption on the smallest angle by a hard-core assumption and then $\Omega$, by the set of admissible configurations $\Omega_{\delta}= \{ \varphi\in\Omega: \forall x,y\in \varphi\times\varphi, x\neq y \quad ||x-y||\geq \delta\}$.

\bibliographystyle{plainnat.bst}
\bibliography{statSpatHAL}

\bigskip \bigskip

\noindent \textbf{Authors}: Jean-Michel Billiot, Jean-Fran\c{c}ois Coeurjolly and Rémy Drouilhet. \\
\noindent \textbf{Address}: LABSAD, BSHM, 1251 avenue centrale BP 47 - 38040 GRENOBLE Cedex 09. \\
\noindent \textbf{E-mail addresses}: 
\begin{itemize}
\item[] \texttt{Jean-Michel.Billiot@upmf-grenoble.fr}
\item[] \texttt{Jean-Francois.Coeurjolly@upmf-grenoble.fr}
\item[] \texttt{Remy.Drouilhet@upmf-grenoble.fr}
\end{itemize}
\noindent \textbf{Corresponding author}: Jean-François Coeujolly.

\newpage

\begin{table}
\begin{tabular}{|l|cc|cc|}
\hline 
Domain $\Lambda_n$& Mean of Estim. of $\theta_2$ &(Std Dev.) & Mean of Estim. of $\theta_3$& (Std Dev.)  \\
\hline
$[-250,250]^2$  & $2.068$ & $0.104$  &$4.382$ & $0.786$\\ 
$[-350,350]^2$  & $2.049$ & $0.071$  &$4.223$ & $0.551$\\ 
$[-450,450]^2$  & $2.041$ & $0.056$  &$4.144$ & $0.436$\\ 
\hline
\end{tabular}
\caption{Empirical mean and standard deviation of maximum pseudo-likelihood estimates of parameters of $\theta_2=2$ and $\theta_3=4$ representing the levels of a multi-Strauss pairwise interaction function where the vector of bounds is assumed to be known and fixed to $\Vect{d}=(0,20,80)$ . These results are obtained from $m=5000$ replications of the point process described by~(\ref{eq-MSDel2}) generated in the domain $[-600,600]^2$. Three sizes of domains $\Lambda_n$ have been considered.} \label{tab-res}
\end{table}

\end{document}